\documentclass{amsart}
\usepackage{amsfonts,latexsym,amssymb,amsmath,amsthm}
\usepackage{enumerate}
\usepackage{graphicx}
\usepackage{float}
\usepackage{caption}

    \theoremstyle{plain}
    \newtheorem{teor}{Theorem}[section]
    
    \newtheorem{coro}[teor]{Corollary}
    \newtheorem{lema}{Lemma}[section]

    \theoremstyle{definition}
    \newtheorem{defin}{Definition}[section]
    \newtheorem{eje}{Example}[section]
    \newtheorem{nota}{Remark}[section]

\numberwithin{equation}{section}

\usepackage{hyperref}
\newcommand{\doi}[1]{\href{http://dx.doi.org/#1}{doi:#1}}

\def\yy{{ y}}

\def\00{{ 0}}

\def\a{ \alpha}

\def\N{\mathbb N}
\def\R{\mathbb R}
\def\C{\mathbb C}
\def\d{\partial}
\def\e{\epsilon}

\usepackage{amsfonts}
\usepackage{amsmath}
\usepackage{xcolor}
\usepackage{bm}

\newcommand{\qbinom}[2]{\genfrac{[}{]}{0pt}{}{#1}{#2}_q}
\newcommand{\qbinomn}[2]{\genfrac{[}{]}{0pt}{}{#1}{#2}_{|q|}}

\begin{document}

\title{$q$-Nagumo norms and formal solutions to singularly perturbed $q$-difference equations}

\author{Sergio A. Carrillo}

\address{Escuela de matem\'{a}ticas, Universidad Nacional de Colombia, Cra. 65 $\#$ 59a-110, Medell\'{i}n, Colombia.}
\email{sacarrillot@unal.edu.co}

\author{Alberto Lastra}

\address{Universidad de Alcal\'{a}, Departamento de F\'{i}sica y Matem\'{a}ticas, Ap. de Correos 20, E-28871 Alcal\'{a} de Henares (Madrid), Spain.}
\email{alberto.lastra@uah.es}

\thanks{Research supported by the project PID2019-105621GB-I00 of Ministerio de Ciencia e Innovaci\'on, Spain. The first author is  partially supported by Universidad Nacional de Colombia under the Proyect 56664(2023). The first author acknowledge Universidad de Alcal\'a (UAH) for its hospitality and support under the ``Giner de los Rios" program (2021) where this work started. The second author is partially supported by Direcci\'on General de Investigaci\'on e Innovaci\'on, Consejer\'ia de Educaci\'on e Investigaci\'on of Comunidad de Madrid, Universidad de Alcal\'a under grant CM/JIN/2021-014, and by Ministerio de Ciencia e Innovaci\'on-Agencia Estatal de Investigaci\'on MCIN/AEI/10.13039/501100011033 and the European Union ``NextGenerationEU''/ PRTR, under grant TED2021-129813A-I00.
}


\maketitle

\begin{abstract}

The aim of this work is to establish the existence, uniqueness and $q$-Gevrey character of formal power series solutions of $q$-analogues of analytic doubly-singular equations. Using a new family of Nagumo norms adapted for $q$-differences we find new types of optimal divergence associated with these problems. We also provide some examples to illustrate our results.
\end{abstract}


\section{Introduction}\label{Sec. Intro}

The genesis of $q$-calculus goes back to Fermat and his computation of  $\int_0^a t^\a dt$, $\a>0$, by subdividing the interval $[0,a]$ using a geometric dissection of ratio $0<q<1$, see \cite[Ch. 10]{andrewaskeyroy}. This led to the development of \textit{Jackson's $q$-integral}, which is the inverse of the $q$-derivative $$d_qf(x):=\frac{f(qx)-f(x)}{qx-x},$$ also known as \textit{Jackson $q$-derivative}. In the current literature, it is common to work with the dilation operator $$\sigma_q(f)(x):=f(qx),$$ instead. Note that both are related by $\sigma_q=(q-1)x d_q+\text{id}$, where $\text{id}$ is the identity. 

% More generally, and loosely speaking this is an example of confluence, meaning that a well-defined object depending on $q$ recovers a usual one when $q\to 1$.

Another origin for $q$-series is Euler's works on partitions leading to $q$-exponential functions and interesting factorizations of these type of series. Other noteworthy $q$-analogues are Gauss' $q$-binomial formula and Heine's $q$-hypergeometric series and their applications in Number Theory, see \cite{Kac}.  These are examples of special functions $y$ satisfying relations of the form $$H(x;y,\sigma_qy,\dots,\sigma_q^ny)=0,$$ known as $q$-difference equations, where $q\in \C^\ast$. Rewriting this equation in terms of $d_q$ and noticing that $d_qf\to f'$ as $q\to 1$, we find that $q$-difference equations can be seen as a discrete counterpart of differential equations. For more historical accounts $q$-calculus and $q$-difference equations  we refer to the survey \cite{DiVizioramisZhang}.

In the analytic setting, linear algebraic $q$-difference equations have been studied since the works of Carmichael \cite{Carmichael}, 
Birkhoff \cite{Birkhoff} and Adams \cite{Adams}. Despite an initial lagged development in comparison with ordinary differential equations (ODEs), Birkhoff's program on the subject has been successfully carried out. This includes the study of symmetries, analytic classification, normal forms, and the inverse Riemann problem. We can mention the works of J.-P. B\'ezivin \cite{bezivin}, L. Di Vizio \cite{vizio}, J. P. Ramis \cite{RamisQ}, J. Sauloy and C. Zhang \cite{ramiszhangSauloy}, and the references therein.

As it is usual in analytic differential problems, divergent formal power series solutions emerge  at irregular singular points. In the $q$-difference framework when $|q|>1$, the divergency is given by a power $s\geq 0$ of  $|q|^{n^2/2}$ and these series are referred as $q$-Gevrey. Optimal values for $s$ are usually found using Newton polygon techniques \cite{canofortuny,marottezhang}, which can be extended to more intricate equations including also differential operators \cite{zhang}. In fact, this is the first step in the study of \textit{summability} of these formal solutions and their Stokes phenomena. Several approaches for $q$-summability have been proposed and developed, taking into account different notions of asymptotic expansions on usual sectors or $q$-spirals. They adapt the use Borel and Laplace transformations by suitable $q$-analogues. For instance, using the Jacobi's theta function \cite{ramiszhang,zhang2}, the two $q$-analogues to the exponential as kernels for the Laplace transform, both with respect to   Riemann's and Jackson's integrals \cite{viziozhang, Dreyfus15, tahara, zhang0}.

Returning to the differential case, after the systematization of summability and its applications in the study of ODEs, the theory was also applied in the setting of singularly perturbed problems, see \cite{CDRSS} and the references therein. However, for problems such as  doubly-singular systems of analytic ODEs of the form \begin{equation}\label{Eq. MS}
    \e^\a x^{p+1} \frac{\d y}{\d x}=F(x,\e,y),
\end{equation} new ideas to identify the correct source of divergence of formal power series solutions were necessary.  In fact, (\ref{Eq. MS}) led to the development of monomial summability in \cite{CDMS}. The key here is to recognize the monomial variable $t=\e^\a x^p$ as the correct one to compute asymptotic expansions. In (\ref{Eq. MS}) we have that  $y=(y_1,\ldots,y_N)\in\C^N$ is a vector of unknown functions, $p,\a\in\N^+$, $F$ is analytic at $(0,0,0)\in\C\times\C\times\C^N$, with $F(0,0,0)=0$, and $DF_y(0,0,0)=0$. The singular perturbation in $\epsilon$ occurs since the nature of the equation changes from differential to implicit as $\e\to 0$. We point out that some $q$-analogues of singularly perturbed ODEs has been studied in \cite{lamasa,malek22} and their references, but only with expansions in the perturbation parameter.

It is precisely (\ref{Eq. MS}) the inspiration of this  work. Our main goal
is to describe the divergence rate of the formal power series solutions to singularly perturbed systems of $q$-difference equations obtained by discretizing (\ref{Eq. MS}). Our aim is aligned with the understanding  of ``\textit{the complete theory of convergence
and divergence of formal series}" \cite[p. 222]{BirkhoffGuenther} for these systems.

More specifically, for $q\in\C$ with $|q|>1$, we consider the problems  
\begin{equation*}
\epsilon^{\alpha}x^{p+1}d_{q,x}(y)(x,\e)=F(x,\e,y),\qquad \epsilon^{\alpha}x^{p}\sigma_{q,x}(y)(x,\e)=F(x,\e,y),
\end{equation*} with similar hypotheses as before. Here we use the notation $\sigma_{q,x}(y)(x,\e):=y(qx,\e)$ and $(q-1)x d_{q,x}+\text{id}=\sigma_{q,x}$  to indicate the action on the first coordinate. However, when the context is clear we will omit this index. 

Although similar, we decided to analyze each one of them by separate, specially because the system involving $d_{q,x}$ is better suited for confluence and it makes sense for $p=-1$. As we mentioned before, although it is common to work only with equations involving $\sigma_q$, a direct approach using $d_q$ can also be fruitful, see, e.g., \cite{essadiq, tahara}. As a matter of fact, $\sigma_q$ and $d_q$ motivate different types of $q$-summability \cite{viziozhang}.

For each equation we will first establish the existence and uniqueness of a solution of the form \begin{equation}\label{Eq. y ini}
\hat{y}(x,\e)=\sum_{n=0}^\infty y_n(\e)x^n=\sum_{n=0}^\infty u_n(x)\e^n.
\end{equation} Then we will determine the growth of the families $\{y_n(\e)\}_{n\in\N}$ and $\{u_n(x)\}_{n\in\N}$ by using majorant series and adequate families of norms, including a new adaptation of Nagumo norms \cite{Nagumo}, that we call \textit{$q$-Nagumo norms}. Incidentally, we can also treat by the same technique both equations in the Fuchsian-like case, namely, when $p=0$. More precisely, we have obtained the following theorem, see below for notation.

\begin{teor}\label{Thm. Main} Fix $q\in \C$ such that $|q|>1$. Consider each one of the systems \begin{align}\label{Eq. principal}
\e^\a x^{p+1}d_{q,x}(y)(x,\e)&=F(x,\e,y),\\
\label{Eq. principalII}
\e^\a x^{p}\sigma_{q,x}(y)(x,\e)&=F(x,\e,y),
\end{align} where $y\in\C^N$, $p\in\N$, $\a\in\N^+$, $F$ is analytic at $(0,0,0)\in\C\times\C\times\C^N$, $F(0,0,0)=0$, and $DF_y(0,0,0)$ is an invertible matrix. Then, each system has a unique formal power series solution $\hat{y}(x,\e)\in \C[[x,\e]]^N$ of the form (\ref{Eq. y ini}). Moreover,

\begin{enumerate}
    \item If $p>0$, there is $r>0$ such that $y_n\in\mathcal{O}_b(D_r)$, $u_n\in\mathcal{O}_b(D_{r/|q|^{\lfloor \frac{n}{\a}\rfloor}})$, and there are constants $C=C(q),A=A(q)>0$ such that, for all $n\geq 0$, $$
\sup_{|\e|\leq r} |y_n(\e)|\leq CA^n |q|^{\frac{n^2}{2p}},\qquad 
\sup_{|x|\leq {r}/{|q|^{\lfloor \frac{n}{\a}\rfloor}} } |u_n(x)|\leq C A^n .$$

    \item If $p=0$, there is $r>0$ such that $y_n\in\mathcal{O}_b(D_{r/|q|^{n/\a}})$, $u_n\in\mathcal{O}_b(D_{r/|q|^{\lfloor n/\a\rfloor}})$, and there are constants $C=C(q),A=A(q)>0$ such that, for all $n\geq 0$, $$
\sup_{|\e|\leq r/|q|^{\frac{n}{\a}}} |y_n(\e)|\leq CA^n,\qquad 
\sup_{|x|\leq {r}/{|q|^{\lfloor \frac{n}{\a}\rfloor}} } |u_n(x)|\leq C A^n.$$
\end{enumerate}
\end{teor}

\

Several remarks are at hand. First, the divergence of the solutions manifests in two ways. Naturally, in some cases the growth of the coefficients is given by a power of the factor $|q|^{n^2/2}$, which is expected for divergence of $q$-difference equations. But also, it is evidenced in the reduction of the radii of the disks where every coefficient is defined. For the coefficients $u_n(x)$ the appearance of the radius $r/|q|^{\lfloor \frac{n}{\a}\rfloor}$ means that, when solved recursively, it is necessary to reduce the radius by a factor for $q$ every $\a$ steps.

The same proof can also be extended to the $q$-difference equation (\ref{Eq. principal}) for $p=-1$ for which we have obtained the following result.

\begin{teor}\label{Prop1}
Fix $q\in \C$ such that $|q|>1$. Consider the system \begin{equation}\label{Eq.3}
\e^\a d_{q,x}(y)(x,\e)=F(x,\e,y),
\end{equation} with the same hypothesis as in Theorem \ref{Thm. Main}. Then, (\ref{Eq.3}) has a unique formal power series solution $\hat{y}=\sum_{n=0}^\infty u_n(x)\e^n$, where $u_n\in\mathcal{O}_b(D_{r/|q|^{\lfloor \frac{n}{\a} \rfloor}})$, for some $r>0$. Moreover,  there are constants $C=C(q),A=A(q)>0$ such that, for all $n\geq 0$, $$\sup_{|x|\leq {r}/{|q|^{\lfloor  
\frac{n}{\a} \rfloor}} } |u_n(x)|\leq C A^n |q|^{n^2/2\a^2}.$$
\end{teor}

\

It is worth mentioning that here we have a higher rate of divergence, in comparison with Theorem \ref{Thm. Main}.
 
These are interesting problems and our results leave an open door to the problem of adapting monomial summability to this context. Several problems are already evident, such as the lose of symmetry in the growth of the coefficients in each variable. In the differential case, when $p=\a=1$, both $y_n$ and $u_n$ are defined in commons disks, for all $n$, and both grow as $n!$. This is no longer valid in the $q$-difference case and new strategies will be necessary to understand how to associate analytic solutions asymptotic to the formal ones obtained here. This is also the reason why we do not expand the solutions in the corresponding monomial or why we do not lift the main equation as in \cite{CaLa23} (technique useful for certain families of holomorphic PDEs). These questions will be addressed in a future work.

The plan for the paper is as follows. In Section \ref{Sec. Preliminaries} we recall the basic facts on some $q$-analogues and the operators $d_q$ and $\sigma_q$, while Section \ref{Sec. Gevrey series} includes some remarks on the spaces of $q$-Gevrey series appearing in our results.  In Section \ref{Sec. Nagumo norms} we introduce  $q$-Nagumo norms, establishing their main properties. Section \ref{Sec. Proofs} is devoted to the proofs of Theorems \ref{Thm. Main} and \ref{Prop1}, where we first describe the general strategy and give details in one case, and only explaining the necessary modifications for the remaining ones. Then, we include several examples in Section \ref{Sec. Examples} to better illustrate our results, showing that the bounds provided by the previous theorems are optimal. Finally, we include an appendix with another possible definition of $q$-Nagumo norms when $q>1$ is real, which are better suited for confluence. In fact, they allow to recover the divergence type of formal solutions of the double singular equation (\ref{Eq. MS}) using confluence, but in general, do not provide optimal bounds, in constrast with the $q$-Nagumo norms of Section \ref{Sec. Nagumo norms}. 

\

\textbf{Notation}. $\N$ is the set of non-negative integers and $\N^+:=\N\setminus\{0\}$. We will be working in $(\C^2,\00)$ with local coordinates $(x_1,x_2)$, $\C[[x_1,x_2]]$ and $\C\{x_1,x_2\}$ will denote the spaces of formal and convergent power series in $(x_1,x_2)$ with  coefficients in $\C$. In the main problem we will also write $(x_1,x_2)=(x,\e)$, to distinguish the singular parameter $\e$. Given $R>0$, we write $D_{R}:=\{x\in\C: |x|<R\}$ for the disc centered at $0\in\C$ and radius $R$. We also set  $\mathcal{O}(\Omega)$ (resp. 
$\mathcal{O}_b(\Omega)$) for the space of $\C$-valued holomorphic (resp. and bounded) functions  on an open domain $\Omega\subseteq\C$. Note that $\mathcal{O}_b(\Omega)$ endowed with the supremum norm is a Banach space.

\section{Preliminaries}\label{Sec. Preliminaries}

In this section we recall some $q$-analogues that will be used in the sequel, including $q$-factorials and $q$-Gevrey series. From now on we fix a complex  number $q\in\C$ with $|q|>1$. We write 
$$[\lambda]_q=\frac{q^\lambda-1}{q-1},\qquad \lambda\in\R,$$ for the $q$-analogue to $\lambda$, which reduces to   $$[n]_q=1+q+\cdots+q^{n-1},\qquad \text{for } n\in\N^+.$$ We have that  $[0]_q=0$ and \begin{equation}\label{Eq. npq}
[np]_q=\frac{q^{p}-1}{q-1} \frac{q^{np}-1}{q^p-1}= [p]_q\cdot [n]_{q^p},\qquad n,p\in\N^+.\end{equation} For future use, we note that $$|[n]_q|\leq [n]_{|q|},\qquad \text{ for } n\in\N,$$ and also that \begin{equation}\label{Eq. nqn asym}
\lim_{n\to+\infty} \frac{[n]_q}{q^n}=\frac{1}{q-1}.
\end{equation}

These constants appear naturally while considering Jackson's $q$-derivative of a function $f$, which is given by $$d_q(f)(x):=\frac{f(qx)-f(x)}{qx-x}=\frac{\sigma_q(f)(x)-f(x)}{qx-x},$$ whenever the expression is defined. As before, $\sigma_q(f)(x):=f(qx)$. For analytic functions $f\in\mathcal{O}(D_r)$, we see that \begin{equation}\label{Eq. Rd}
    \sigma_q(f), d_q(f)\in \mathcal{O}(D_{r/|q|})
\end{equation} and they can be computed term by term using its power series expansion according to the rules $$d_q(x^n)=[n]_q x^{n-1},\qquad \sigma_q(x^n)=q^n x^n,\qquad n\in\N.$$ On the other hand, this formula allows to consider $d_q, \sigma_d:\C[[x]]\to\C[[x]]$, also defined term by term. In this setting, Leibniz rule is replaced by \begin{equation}\label{qLeib}d_q(fg)(x)=d_q(f)(x)g(x)+f(qx)d_q(g)(x).\end{equation}

We recall the coefficients 
\begin{equation}\label{Eq. a,q}
(a;q)_n=\prod_{j=0}^{n-1} (1-aq^j),\qquad (a;q^{-1})_\infty=\prod_{j=0}^{\infty} (1-aq^{-j}),\qquad a\in\C.
\end{equation} The second one is convergent as we can compare it with a geometric series. The $q$-factorial is defined accordingly as $$[n]_q^!=[1]_q [2]_q\cdots [n]_q=\frac{(q;q)_n}{(1-q)^n}.$$ In general, for $|q|>1$, since  $\lambda\in\R\longmapsto [\lambda]_{|q|}$ is a strictly increasing function, the same holds for the map $n\in\N\longmapsto [n]^!_{|q|}$. Therefore,  \begin{align*}
\frac{[n-p]_{|q|}}{([n]^!_{|q|})^{1/p}}  = \frac{[n-p]_{|q|}}{([n]_{|q|}\cdots[n-p+1]_{|q|})^{1/p} ([n-p]^!_{|q|})^{1/p}}&\leq \frac{[n-p]_{|q|}}{[n-p+1]_{|q|} } \frac{1}{([n-p]^!_{|q|})^{1/p}},
\end{align*} and thus 
\begin{equation}\label{Eq. Inq n-p qfactor}
\frac{[n-p]_{|q|}}{([n]^!_{|q|})^{1/p}}    \leq \frac{1}{([n-p]^!_{|q|})^{1/p}},\qquad \text{for integers } n>p>0.\end{equation}

Another useful $q$-analogue is the $q$-binomial coefficient $$\qbinom{n}{j}
= \frac{[n]^{!}_{q}}{[j]^{!}_{q} [n-j]^{!}_{q}},\qquad 0\leq j\leq n,$$ which is a polynomial in $q$ of degree $j(n-j)$ and satisfies the $q$-analogue to Pascal's formula 
\begin{equation}\label{Eq. q binom}
    \qbinom{n}{j}=\qbinom{n-1}{j-1}+q^j \qbinom{n-1}{j},
\end{equation} see, e.g., \cite[Chapter 6]{Kac}. In particular, it follows that $\qbinomn{n}{j}\geq 1$, i.e., \begin{equation}\label{Eq. Inq qbinom}
[j]^!_{|q|} [n-j]^!_{|q|}\leq [n]^!_{|q|}.
\end{equation} which will be used later.

\begin{nota}\label{Rmk. logcovex} In general,  a sequence $\{M_n\}_{n\geq 1}$ of positive real numbers such that $M_0=1$, satisfies the inequality $M_nM_k\leq M_{n+k}$, for all $n,k\geq 0$, when it is \textit{log-convex}, i.e., when it holds that $$M_n^2\leq M_{n-1} M_{n+1},$$ for all $n\geq 1$. This fact follows easily using induction on $k$. An example of this situation is precisely  the sequence $M_n=([n]_{|q|}^!)^s$, where $s>0$ and $|q|>1$. Indeed, $M_{n}/M_{n-1}=[n]_{|q|}^s$ which is increasing in $n$.
\end{nota}

In general, $q$-factorials determine the  divergence rate of solutions of singular $q$-difference equations. Two classical examples are the following.

\begin{eje} A $q$-analogue to Euler's equation is $$x^2d_qy(x)+y(x)=x,$$ having as unique  formal power series solution $$\hat{E}_q(x):=\sum_{n=0}^\infty (-1)^n [n]_q^! x^{n+1},$$ which is divergent for $|q|>1$. In contrast, the unique formal power series solution of the $q$-difference equation $$x\sigma_qy(x)=y(x)-1$$ is given by  $$\hat{Y}_q(x)=\sum_{n=0}^\infty q^{n(n+1)/2} x^{n+1},$$ which also diverges for $|q|>1$. As $|q|\to 1^+$, the divergence of $\hat{E}(x)=\sum_{n=0}^\infty n! x^{n+1}=\lim_{|q|\to 1^+} \hat{E}_q(x)$ persists, while $\lim_{|q|\to 1^+} \hat{Y}_q(x)=\frac{x}{1-x}$ converges.
\end{eje}

In this framework, the following notion plays the role of Gevrey series.

\begin{defin} A series $\hat{f}=\sum_{n\geq 0} a_n x^n\in \C[[x]]$ is \textit{$s$-$q$-Gevrey}, where $s\geq 0$ and $|q|>1$, if there are constants $C=C(q),A=A(q)>0$ such that \begin{equation}\label{Eq. Def s q z}
|a_n|\leq CA^n |q|^{sn^2/2}.
\end{equation} We will denote the space of such series by $\C[[x]]_{q,s}$.
\end{defin}

We remark that the $q$-factorial is desirable while working with $q$-difference equations involving $d_q$ and the coefficient $|q|^{n^2/2}$ is natural when the equation is written in terms of $\sigma_q$.

\begin{nota} We can interchange the term $|q|^{n^2/2}$ above for $|q|^{n(n\pm 1)/2}$. Also, we can use the terms $[n]_{|q|}^!$ or $|[n|_q^!|$, i.e., asking for a series to be $s$-$q$-Gevrey is equivalent to  request that $$|a_n|\leq D B^n ([n]_{|q|}^!)^s,\qquad n\in\N,$$ for some constants $B,D>0$ that might also depend on $q$. Indeed, since  $$[n]_q=q^{n-1}[n]_{q^{-1}}, \quad \text{ and }\quad [n]_{q}^{!}=q^{n(n-1)/2} [n]_{q^{-1}}^{!},$$ we have the limit $$\lim_{n\to+\infty} [n]_{q}^{!}/{\displaystyle \frac{q^{\frac{n(n+1)}{2}}}{(q-1)^n}}=(q^{-1},q^{-1})_\infty=c(q).$$ This means that $[n]_{q}^{!}$ is asymptotically equivalent to $ c(q)(q-1)^{-n}q^{n(n+1)/2}$, c.f., \cite[p. 55]{RamisQ}. Thus, $[n]_{|q|}^!$ is asymptotically equivalent to $ c(|q|)(|q|-1)^{-n} |q|^{n(n+1)/2}$, and $|[n|_q^!|$ is to $ |c(q)| |q-1|^{-n} |q|^{n(n+1)/2}$. We also note that we have the simple explicit bound $$\frac{c(|q|)}{(|q|-1)^{n}}|q|^{\frac{n(n+1)}{2}}\leq [n]_{|q|}^!\leq \frac{1}{(|q|-1)^n} |q|^{\frac{n(n+1)}{2}},$$  relating $|q|^{\frac{n(n+1)}{2}}$ and $[n]_{|q|}^!$, up to geometric terms.
\end{nota}

It is easy to check that $\C[[x]]_{q,s}$ is stable under sums, products and the operators $d_q$ and $\sigma_q$. Regarding ramifications and changes in $s$, we highlight that:
\begin{enumerate}
    \item $\hat{f}$ is $s$-$q$-Gevrey if and only if it is $1$-$q^s$-Gevrey.
    
    \item Fix $p\in\N^+$ and write $\hat{f}(x)=\sum_{j=0}^{p-1} x^j \hat{f}_j(x^p)$, where $\hat{f}_j(t)=\sum_{n=0}^\infty a_{np+j} t^n$. Then $\hat{f}(x)$ is $s$-$q$-Gevrey (in $x$) if and only if each $\hat{f}_j(t)$ is $sp$-$q^p$-Gevrey (in $t$). In fact, we obtain from (\ref{Eq. Def s q z}) that $$|a_{np+j}|\leq CA^{np+j} |q|^{s(np+j)^2/2} = C_j A_j^n |q|^{sp^2n^2/2},$$ where $C_j=CA^j|q|^{sj^2/2}$ and $A_j=A^p |q|^{pjs}$. The converse can also be easily proved.
\end{enumerate}

\section{Some spaces of \texorpdfstring{$q$}--Gevrey series}\label{Sec. Gevrey series}

This section introduces two spaces of $q$-Gevrey series in two variables describing the phenomena encountered in Theorem \ref{Thm. Main}. 

Let $\hat{f}\in \C[[x_1,x_2]]$ be a formal power series,  written canonically as \begin{equation}\label{Eq. f ex1}
\hat{f}=\sum_{n,m=0}^\infty a_{n,m} x_1^n x_2^m = \sum_{n=0}^\infty y_n(x_2) x_1^n=\sum_{m=0}^\infty u_m(x_1) x_2^m,    
\end{equation} for $y_n\in\C[[x_1]]$ and $u_m\in\C[[x_2]]$. 

% Moreover, $\hat{f}$ can also be expanded as $$\hat{f}=\sum_{n=0}^\infty f_n(x_1,x_2) (x_1 x_2)^n,\qquad f_n=a_{n,n}+\sum_{m=1}^\infty  a_{n+m,n} x_1^m +a_{n,n+m}x_2^m,$$ writing it as a power series in the monomial $t=x_1x_2$. 

First, we consider the space $$\mathcal{O}_{0}^q:=\left\{\hat{f}\in\C[[x_1,x_2]] : 
 |a_{n,m}|\leq CA^{n+m} |q|^{nm}, \text{ for some } C,A>0 \text{ and all }n,m \right\},$$ that in some sense plays the roles of the ring $\mathcal{O}_0=\C\{x_1,x_2\}$ of convergent power series, for the case $q=1$. We can characterize their elements as follows.

 \begin{lema}\label{Lm. 1} Let $\hat{f}\in\C[[x_1,x_2]]$ as in (\ref{Eq. f ex1}). The following are equivalent:\begin{enumerate}
     \item $\hat{f}\in \mathcal{O}_{0}^q$,

     \item For some $r,B,D>0$, $y_n\in \mathcal{O}_b(D_{r/|q|^n})$, for all $n$, and \begin{equation*}\label{Eq. f O0qA}
\sup_{|x_2|\leq r/|q|^n} |y_n(x_2)|\leq DB^n,\end{equation*}

     \item For some $r,B,D>0$, $u_m\in \mathcal{O}_b(D_{r/|q|^m})$, for all $m$, and \begin{equation*}\label{Eq. f O0qB}
\sup_{|x_2|\leq r/|q|^m} |u_m(x_1)|\leq DB^m.\end{equation*} \end{enumerate}  
 \end{lema}

\begin{proof} To show that (1) implies (2) and (3), let  $\hat{f}\in \mathcal{O}_{0}^q$. Then $$|y_n(x_2)|\leq \sum_{m=0}^\infty CA^{n+m}  |q^nx_2|^{m}\leq 2CA^n,$$ if $|x_2|\leq A/2|q|^n$, and analogously for $u_m(x_1)$. Now, if (2) holds, by Cauchy's inequalities, $$|a_{n,m}|=\left|\frac{1}{m!}\frac{\d^m y_n}{\d x_2^m}(0)\right|\leq \frac{DB^n}{(\rho/|q|^n)^m}=DB^n\rho^{-m} |q|^{nm},\qquad \text{ for } 0<\rho<r,$$ and for all $n,m\in\N$.  Thus (1) is valid. In the same way, (3) implies (1). Alternatively, note the symmetry in the previous definition: $\hat{f}(x_1,x_2)\in \mathcal{O}_{0}^q$ if and only if $\hat{f}(x_2,x_1)\in \mathcal{O}_{0}^q$, so we can exchange the role of the variables.
\end{proof}

Thanks to Lemma \ref{Lm. 1} we see that Theorem \ref{Thm. Main} (2) claims that the unique solutions of the corresponding systems belongs to $\mathcal{O}_0^{q^{1/\a}}$.

\begin{eje}\label{Ex. Model} Consider the series $$M(x_1,x_2):=\sum_{n,m=0}^\infty q^{nm} x_1^n x_2^m=\sum_{n=0}^\infty \frac{x_1^n}{1-q^nx_2}=\sum_{m=0}^\infty \frac{x_2^m}{1-q^m x_1}\in\mathcal{O}_0^q,$$ which is confluent to the geometric series $\left[(1-x_1)(1-x_2)\right]^{-1}$ when $q\to 1$. Although $M$ does not define an analytic function at $(0,0)$, given any $\delta>0$ the series converges on the set
$$U_{q,\delta}^1:=U_{q,\delta}\times D_{\delta},\qquad U_{q,\delta}:=\left\{x_1\in\C\,:\, |1-q^mx_1|>\delta^m\hbox{ for all } m\ge0\right\}.$$
The second variable consists of a disc of radius $\delta>0$. In the first variable we have the intersection of the complement of circles with radius $(\delta/|q|)^m$ and center $1/|q|^m$. If $\delta=1$, these circles have $0$ at the boundary. If additionally $q$ is real,  then $U_{q,1}=\{x_1\in\C : |1-x_1|>1\}$ since $D(1/q^{m+1},1/q^{m+1})\subset D(1/q^m,1/q^m)$. For general values of $q$ and $\delta$, these circles spiral around the origin, see Figure \ref{Fig.Spirals} for an example. Since $M(x_1,x_2)=M(x_2,x_1)$, $M$ also defines an analytic function on $U_{q,\delta}^2=D_\delta\times U_{q,\delta}$, by interchanging the order of the variables. 

Finally, we note that $M(x,\e)$ is the unique formal power series solution of $$\e \sigma_{q,x}y=y-\frac{1}{1-x},$$ corresponding to equation (\ref{Eq. principalII}) in Theorem \ref{Thm. Main}(2).

\begin{center}
\begin{figure}
\includegraphics[scale=0.22]{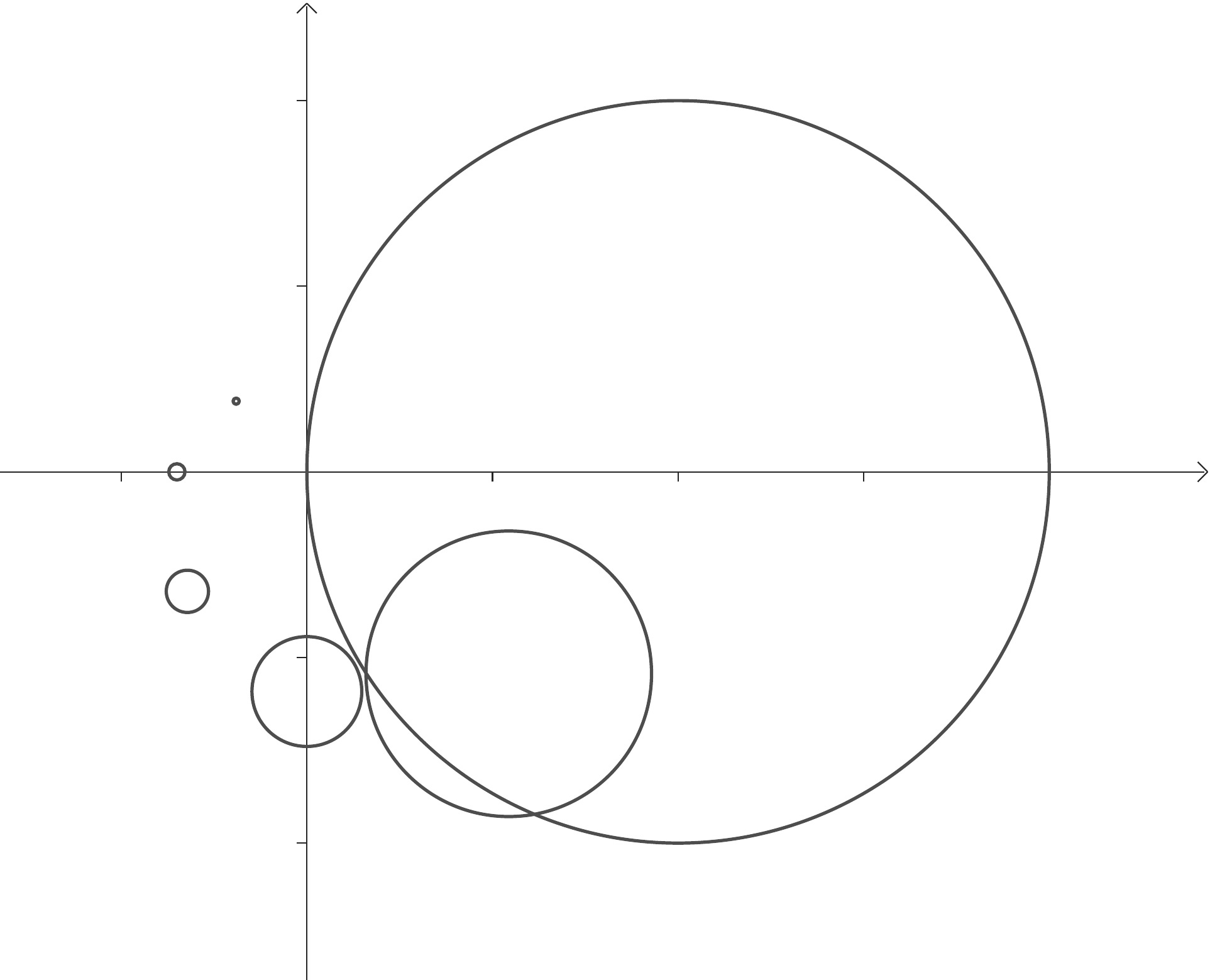}	%\vspace{-0.9cm}	
	
\caption{Circles in $U_{q,\frac{1}{2}}$, for $q=\frac{3}{2}e^{i\pi/4}$.}
	\label{Fig.Spirals}
\end{figure}
\end{center}

\end{eje}

The second type of series that appear naturally in Theorem \ref{Thm. Main} are the following, named after the analogy with the differential case.

\begin{defin}\label{Def. q mon series} Fix a monomial  $x_1^px_2^\a$. We say that a series $\hat{f}=\sum a_{n,m} x_1^n x_2^m\in\C[[x_1,x_2]]$ is \textit{$1$-$q$-Gevrey in $x_1^px_2^\a$} if there are constants $C=C(q),A=A(q)>0$ such that $$|a_{n,m}|\leq CA^{n+m}\min\{ |q|^{n^2/2p},  |q|^{nm/\a} \}.$$  The space of these series will be denoted by $\mathcal{O}^q_{x_1^p x_2^\a}$.
\end{defin}

Proceeding as in Lemma \ref{Lm. 1}, we can establish the following characterization for elements of $\mathcal{O}^q_{x_1^p x_2^\a}$.

\begin{lema}\label{Lm. 2} Let $\hat{f}\in\C[[x_1,x_2]]$ as in (\ref{Eq. f ex1}). The following assertions are equivalent:\begin{enumerate}
     \item $\hat{f}\in \mathcal{O}^q_{x_1^p x_2^\a}$,

     \item For some $r,B,D>0$, $y_n\in \mathcal{O}_b(D_{r})$, $u_n\in \mathcal{O}_b(D_{r/|q|^{n/\a}})$ for all $n$, and \begin{equation*}\label{Eq. f O0qC}
\sup_{|x_2|\leq r} |y_n(x_2)|\leq DB^n|q|^{n^2/2p},\qquad \sup_{|x_1|\leq r/|q|^{n/\a}} |u_n(x_1)|\leq DB^n.\end{equation*}
\end{enumerate}
 \end{lema}

We collects below some algebraic properties on these spaces including the behavior under ramifications. The proof is straightforward.

\begin{lema} Consider $|q|>1$ and a monomial $x_1^p x_2^\a$. The following assertions hold:\begin{enumerate}
\item $\mathcal{O}^{q}_{0}$ and $\mathcal{O}^q_{x_1^p x_2^\a}$ are rings stable by the operators $d_{q,x_j}, \sigma_{q,x_j}$, $j=1,2$

\item $\mathcal{O}^q_{x_1^p x_2^\a}\subset \mathcal{O}^{q^{1/\a}}_{0}$. 

\item $\hat{f}=\sum_{j=0}^{p-1} \hat{f}_j(x_1^p,x_2) x_1^j\in \mathcal{O}_{0}^q$ if and only if $\hat{f}_j(z,x_2)\in \mathcal{O}_{0}^{q^{p}}$, for $j=0,1,\dots,p-1$.

\item $\hat{f}=\sum_{0\leq j<p, 0\leq l<\a} \hat{f}_{j,l}(x_1^p,x_2^\a) x_1^j x_2^l\in \mathcal{O}^q_{x_1^px_2^\a}$ if and only if $\hat{f}_j(z,\eta)\in \mathcal{O}_{z\eta}^{q^{p}}$, for all $j=0,1,\dots,p-1$, $l=0,1,\dots,\a-1$. 
\end{enumerate}
    
\end{lema}

We conclude this section remarking that the series $\hat{y}=\sum a_{n,m} x_1^n x_2^m$ appearing in Theorem \ref{Prop1} are precisely those satisfying bounds of the type $$|a_{n,m}|\leq CA^{n+m} |q|^{nm/\a}\cdot |q|^{n^2/2\a^2},$$ for certain constants $C=C(q), A=A(q)>0$. In contrast, this is a higher divergence rate compared to the one describes in the spaces $\mathcal{O}_o^q$ and $\mathcal{O}_{x_1^px_2^\a}^q$.

% $$|b_n(x_1)|\leq CA^{2n} q^{n^2/2}n! \sum_{m=0}^\infty A^mq^{n(n+m)} x_1^m=\frac{CA^{2n}}{1-Aq^n|x_1|} q^{3n^2/2}n!$$ $$|c_n(x_2)|\leq  \frac{CA^{2n}}{1-A|x_2|} [n]!_q$$

% On the other hand, we introduce the space $$\mathcal{O}_{r,q}:=\{\hat{f}\in \C[[x_1,x_2]] : y_n\in\mathcal{O}_b(D_r), u_n\in\mathcal{O}_b(D_{r/|q|^n}), \text{ for all }n\},$$ for some fixed $r>0$. Note that for $q=1$ this is simply saying that the coefficients in (\ref{Eq. f ex1}) have a common radius of convergence. In our case, the radius in the first variable will reduce geometrically with $q$ due to the action of $d_{q,x_1}$ or $\sigma_{q,x_1}$. In fact, if $\hat{f}\in\mathcal{O}_{r,q}$, then $x_2\cdot d_{q,x_1}(\hat{f}), x_2\cdot  \sigma_{q,x_1}(\hat{f})\in\mathcal{O}_{r,q}$, thanks to (\ref{Eq. Rd}).

\section{The \texorpdfstring{$q$}--Nagumo norms}\label{Sec. Nagumo norms}

In order to treat the divergence phenomenon generated by the singular parameter in equations (\ref{Eq. principal}), (\ref{Eq. principalII}), and (\ref{Eq.3}),  we introduce and develop an adaptation of modified Nagumo norms. These were introduced in \cite{CDRSS} to obtain the divergence rate of solutions of singular perturbations of ODEs. In turn, the former were based on  classical Nagumo norms, first introduced in M. Nagumo in his work \cite{Nagumo} on power series solutions of analytic PDEs.

\begin{figure}

\centering{
\includegraphics[scale=0.15]{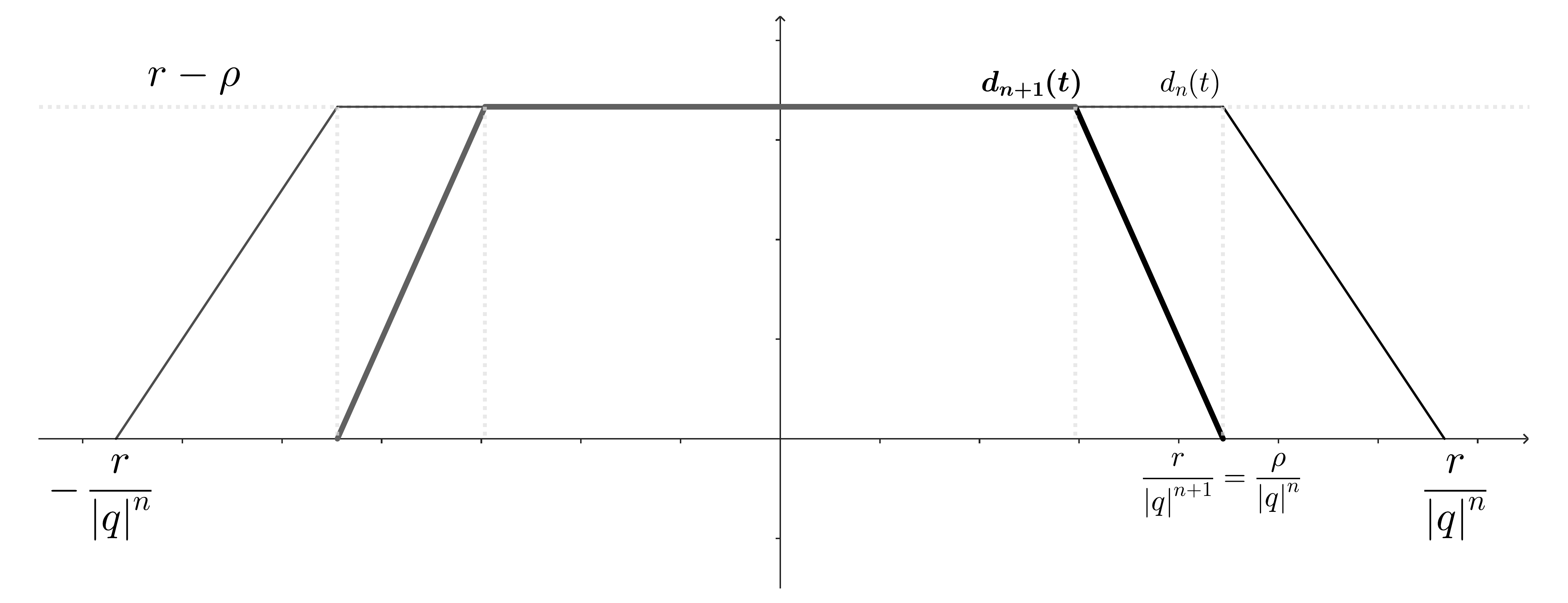}	}
	
	\caption{The auxiliary functions $d_n(t)$.}
	\label{Fig.Dn}
\end{figure}

Fix $|q|>1$, $0<\rho<r$, and $n\in\N$. Consider $d_n:D_{r/|q|^n}\to\R$ given by $$d_n(t)=\begin{cases}
r-|q^nt|, & \rho/|q|^n\leq |t|\leq r/|q|^n,\\
r-\rho, & |t|\leq \rho/|q|^n.
\end{cases}$$ Choosing $\rho$ adequately, the following assertions hold, see Figure \ref{Fig.Dn}.

\begin{lema}\label{Lema dn} If $\rho=r/|q|$, then \begin{enumerate}
    \item $|d_n(t)-d_n(s)|\leq |q|^n|t-s|$.

    \item $d_n(qt)=d_{n+1}(t)$, for $|t|\leq r/|q|^{n+1}$.
    
    \item $d_{n+1}(t)\leq d_n(t)\leq r(1-1/|q|)$, for all $n\in\N$ and $|t|\leq r/|q|^{n+1}$.
\end{enumerate}
\end{lema}

\noindent From now on we fix the value $\rho=r/|q|$ and consider $d_n$ accordingly. For $f\in\mathcal{O}_b(D_r)$, we define the $n$th \textit{$q$-Nagumo norm} by  $$\|f\|_n:=\sup_{|x|\leq r/|q|^n} |f(x)|\cdot  d_n(|x|)^n.$$ To simplify notation we omit the dependence on $r$ and $q$. Note that for $n=0$, $\|f\|_0$ is simply the supremum norm. Also, if $\|f\|_n$ is finite, we have that $$|f(x)|\leq \frac{\|f\|_n}{d_n(|x|)^n},\qquad \text{ for } |x|<r/|q|^n.$$ This family of norms takes into account that factor $r-q^n|x|$ that measures correctly the distance of $x$ to the boundary of the disk $D_{r/|q|^n}$. In this way, we obtain information on the operators $d_q$ and $\sigma_q$, as we shall see it in the following lemma.

\begin{lema} Let $|q|>1$, $n,m\in\N$ and $f,g \in\mathcal{O}_b(D_r)$. The following assertions hold:\begin{enumerate}
    \item $\|f+g\|_n\leq \|f\|_n+\|g\|_n$ and $\|fg\|_{n+m}\leq \|f\|_n \|g\|_m$.
    
    \item $\|d_q(f)\|_{n+1}\leq 2|q|^{n+1}\|f\|_n$.
    
    \item $\|\sigma_q(f)\|_{n+1}\leq r(1-\frac{1}{|q|})\|f\|_n$.
\end{enumerate}
\end{lema}

\begin{proof} For (1), using Lemma \ref{Lema dn}(3), we see that 
$$|f(x)g(x)|d_{n+m}(|x|)^{n+m}\leq |f(x)|d_n(|x|)^n\, |g(x)| d_m(|x|)^m\leq\|f\|_n \|g\|_m,$$ for all $|x|\leq r/|q|^{n+m}$. This proves the inequality for the product. 

For (2) fix $x$ such that $|x|\leq r/|q|^{n+1}$. To bound $d_q(f)$, consider first the case $|x|\geq \rho/|q|^{n+1}=r/|q|^{n+2}$. Then \begin{equation}\label{Eq. aux Lema main}\left|\frac{f(qx)-f(x)}{x}\right|\leq \frac{|q|^{n+2}}{r}\|f\|_n\left(\frac{1}{d_n(|qx|)^n}+\frac{1}{d_n(|x|)^n}\right)\leq \frac{2 |q|^{n+2} \|f\|_n}{r d_{n+1}(|x|)^n},\end{equation} thanks to Lemma \ref{Lema dn}(2) and (3). Thus, we find that \begin{align*}
    \left|\frac{f(qx)-f(x)}{x}\right| d_{n+1}(|x|)^{n+1}
    &\leq 2\frac{|q|^{n+2}}{r}\|f\|_n d_{n+1}(|x|)\\
    &\leq 2|q|^{n+2}\|f\|_n \frac{(r-\rho)}{r}=2|q|^{n+1}\|f\|_n(|q|-1).
\end{align*}

In the case $|x|\leq\rho/|q|^{n+1}$, using the Maximum Modulus Principle and (\ref{Eq. aux Lema main}) we see that \begin{align*}
    \left|\frac{f(qx)-f(x)}{x}\right| &\leq \sup_{|x|=\rho/|q|^{n+1}} \left|\frac{f(qx)-f(x)}{x}\right|\leq \frac{2 |q|^{n+2} \|f\|_n}{r (r-\rho)^n},
    \end{align*} and therefore, \begin{align*}
    \left|\frac{f(qx)-f(x)}{x}\right| d_{n+1}(|x|)^{n+1}
    &\leq 2|q|^{n+1}\|f\|_n(|q|-1).
\end{align*} Since $|q|-1\leq |q-1|$, these inequalities establish (2).

Finally, to prove (3) simply note that by definition of $\|f\|_n$ we have $$|f(qx)|\leq \frac{\|f\|_n}{d_n(|qx|)^{n}}=\frac{\|f\|_n}{d_{n+1}(|x|)^{n}},\qquad \text{ for } |x|<r/|q|^{n+1}.$$ In this way, $|f(qx)| d_{n+1}(|x|)^{n+1}\leq \|f\|_n d_{n+1}(|x|)<\|f\|_n(r-\rho)$ as required.
\end{proof}

\begin{nota}\label{Rmk Nagumo vm} The previous lemma can be extended to vector- and matrix-valued analytic maps in a straightforward way. Indeed, if $f\in\mathcal{O}_b(D_r)$,  $z=(z_1,\dots,z_N)\in\mathcal{O}_b(D_r)^N$, and $A=(a_{ij})\in \mathcal{O}_b(D_r)^{N\times N}$, setting $$\|z\|_n=\max_{1\leq j\leq N} \|z_j\|_n,\qquad \|A\|_n=\max_{1\leq i\leq N} \sum_{j=1}^N \|a_{i,j}\|_n,$$ it follows that $$\|f\cdot z\|_{n+m}\leq \|f\|_n \|z\|_m, \quad\|A\cdot z\|_{n+m}\leq \|A\|_n \|z\|_m,\quad \|d_qz\|_{n+1}\leq 2|q|^{n+1} \|z\|_n,$$ for all $n,m\geq 0$.
\end{nota}

\section{Proof of the main results}\label{Sec. Proofs}

This section is devoted to prove Theorems \ref{Thm. Main} and \ref{Prop1} regarding the divergence rate of the formal power series solutions of equations (\ref{Eq. principal}), (\ref{Eq. principalII}), and (\ref{Eq.3}), respectively. The proof follows the usual majorant series technique. For  expansions in the parameter $\epsilon$, we use the $q$-Nagumo norms to control the action of $d_{q,x}$ and $\sigma_{q,x}$ that emerges in the recurrences induced by the main equations.

The structure of the proofs is the same in all cases and it has been successfully used in the differential setting, see e.g., \cite{CarrVF, cahu, zhangNagumo}. In order to prove the results we will apply the following steps:

\begin{enumerate}
    \item Apply rank reduction in $\epsilon$ to assume $\a=1$.
    
    \item Find the formal solution $\hat{y}$, starting by searching the first coefficient (as a function of $x$ or $\epsilon$) using the Implicit Function Theorem. Then, set the recurrences to  determine the other terms.
    
    \item Using the supreme norm in $D_r$ or $D_{r/|q|^n}$ ($\hat{y}$ as power series in $x$) or the $q$-Nagumo norms ($\hat{y}$ as power series in $\epsilon$), find a system of inequalities for the norms $z_n\geq 0$ of the coefficients. Then, divide them by an adequate log-convex sequence $\{M_n\}_{n\geq 1}$ to find a sequence of inequalities for $z_n/M_n$.
    
    \item Finally, associate an analytic problem having a unique convergent power series solution $\hat{w}=\sum w_n \tau^n$ that majorises the series $\sum \frac{z_n}{M_n}\tau^n$, i.e., $z_n/M_n\leq w_n$, for all $n$.
    \end{enumerate}

Except for Step $3$ where $M_n$ is chosen, the arguments are quite similar in all cases. Therefore, we will only write one in detail, namely, determining the growth of the coefficients of $\hat{y}$ as a power series in $x$. For the others, we will only indicate the necessary modifications.

\begin{proof}[Proof of Theorem \ref{Thm. Main}]

\textit{Solution as a power series in $x$. Case $p>0$}.

\textbf{Step 1: Rank reduction}. We can perform rank reduction on the systems (\ref{Eq. principal}) and (\ref{Eq. principalII}), in exactly the same way as for ODEs, to reduce this to the case $\a=p=1$. However, to preserve the nature of the problem we only reduce rank in $\epsilon$, see Remark \ref{Rmk RR}. Thus, we write \begin{equation}\label{Eq. y RRe}
y(x,\e)=\sum_{j=0}^{\a-1} y_j(x,\epsilon^\a) \e^j,\qquad w(x,\eta)=(y_0(x,\eta),y_1(x,\eta),\dots,y_{\a-1}(x,\eta)).\end{equation} Notice that $\e^\a x^{p+1}d_{q,x}y(x,\e)=\sum_{j=0}^{\a-1} \eta x^{p+1}d_{q,x}y_j(x,\eta) \e^j$ and $\e^\a x^{p}\sigma_{q,x}y(x,\e)=\sum_{j=0}^{\a-1} \eta x^{p}\sigma_{q,x}y_j(x,\eta) \e^j$. Also, expanding $$F(x,\e,y)=\sum_{l=0}^{\a-1} F_l\left(x,\e^\a,\sum_{k=0}^{\a-1} y_k(x,\epsilon^\a) \e^k\right) \e^l = \sum_{j=0}^{\a-1} \widetilde{F}_j(x,\eta,w) \e^j,$$ and setting $G(x,\eta,w)=(\widetilde{F}_0, ,\dots,\widetilde{F}_{\a-1})$, it follows that $D_wG(0,0,0)$ is a block-diagonal matrix of size $N\a$ having all diagonal blocks equal to $D_yF(0,0,0)$, see \cite[Section 4]{CDMS} for details in the differential case. Thus, this matrix is invertible. The new systems, of dimension $N\a$, take the form $$\eta x^{p+1} d_qw=G(x,\eta,w),\qquad \eta x^{p} \sigma_{q}w=G(x,\eta,w)$$ and have the same structure as (\ref{Eq. principal}) and (\ref{Eq. principalII}), but now $\a=1$.

\textbf{Step 2: Existence}. We can assume now that  (\ref{Eq. principal}) and (\ref{Eq. principalII}) have the form 
\begin{equation}\label{Eq. principal2}
\e x^{p+1}d_{q}(y)(x,\e)=F(x,\e,y),\qquad 
 \e x^{p}\sigma_{q}(y)(x,\e)=F(x,\e,y)\end{equation} where $y\in\C^N$, $F$ is holomorphic near $(0,0,0)\in \C\times\C\times\C^N$,  $F(0,0,0)=0$, and $DF_y(0,0,0)$ is invertible.

The existence and uniqueness of a solution $\hat{y}\in\C[[x,\e]]^N$ follows by simply replacing the power series $\hat{y}$ into (\ref{Eq. principal2}) and solving the coefficients recursively. This can be done thanks to the hypothesis of $DF_y(0,0,0)$ being invertible. Note that $\hat{y}(0,0)=0$, since $F(0,0,0)=0$. Here we are interested in the coefficients $y_n(\epsilon)$ of $\hat{y}$ when written in the first form of (\ref{Eq. y ini}). To establish the recurrences that define them, let us write $${F}(x,\e,{y})=b(x,\e)+A(x,\e){y}+H(x,\e,y),\qquad H(x,\e,y)=\sum_{I\in\N^N, |I|\geq 2} A_I(x,\e) {y}^I,$$ as a convergent power series in $\yy$ with coefficients $A(x,\e)=\sum_{n=0}^\infty A_{n*}(\e)x^n$, and analogously for $b(x,\e)$ and each $A_I(x,\e)$. Thus, there is $r>0$ such that $A_{n*}\in \mathcal{O}_b(D_r)^{N\times N}$ and $b_{n*}, A_{I,n*}\in \mathcal{O}_b(D_r)^N$, for all $n\in\N$. Note that we are using the notation $y^{I}=y_1^{i_1}\cdots y_N^{i_N}$, where $I=(i_1,\dots,i_N)\in\N^N$.

The first coefficient $y_0(\e)$ is determined by solving the implicit equation $$F(0,\e,y_0(\e))=0.$$ This has a unique analytic solution $y_0(\e)\in \mathcal{O}_b(D_{r'})^N$, for some $r'>0$, via the Implicit Function Theorem since  $DF_y(0,0,0)=A(0,0)$ is invertible. Reducing $r$ if necessary, we can assume that $r'=r$ and that $A(0,\e)=A_{0*}(\e)$ is also invertible and bounded for all $|\e|<r$.

After the the change of variables $y\mapsto y-y_0(\e)$ in the initial equation (\ref{Eq. principal2}), we can assume now that $y_0(\e)=0$ Thus we obtain a similar $q$-difference equation such that $F(0,\e,0)=b_{0*}(\e)=0$. Now, after replacing $\hat{y}=(\hat{y}_1,\dots,\hat{y}_N)=\sum_{n=1}^\infty y_n(\e) x^n$, where $y_n(\e)=( y_{1,n}(\e),\dots,y_{N,n}(\e))$, into (\ref{Eq. principal2}), we obtain the recurrence \begin{align}\label{Eq. principal en x}
\e [n-p]_q y_{n-p}(\e)=&b_{n*}(\e)+\sum_{j=1}^n A_{n-j*}(\e) y_{j}(\e)\\
\nonumber &+\sum_{k=2}^n \sum_{\ast_k} \prod_{{1\leq l\leq N}\atop{1\leq j\leq i_l}} y_{l,n_{l,j}}(\e)\,\cdot  A_{I,m*}(\e),\quad n\geq 1,
\end{align} for the first system. The same recurrence holds for the second system with $[n-p]_{|q|}$ replaced by $|q|^{n-p}$. The inner sum indicated with $(\ast_k)$ is taken over all $I=(i_1,\dots,i_N)\in\N^N$ such that $|I|=k$, $m$ satisfying $0\leq m\leq n-k$, and $n_{l,j}\geq 1$ such that $n_{1,1}+\cdots+n_{1,i_1}+\cdots+n_{N,1}+\cdots+n_{N,i_N}+m=n$. Note in particular that $n_{l,j} < n-m\leq n$, thus no component of $\hat{y}_n$ appears in the former sum. Finally, the left-side of the equation is understood as zero for $n<p$.

Since $A_{0*}(\e)$ is invertible for $|\e|<r$, the coefficient $y_{n}(\e)$ is uniquely determined recursively by (\ref{Eq. principal en x}) and is analytic and bounded on $D_r$. The same holds for the system with the operator $\sigma_{q}$.

\textbf{Step 3: Estimates for $y_n(\e)$}. We use the supreme norm $\|\cdot\|$, i.e., the $q$-Nagumo norms of order $0$ as in Remark \ref{Rmk Nagumo vm}. Let $c=\|A_{0*}^{-1}\|>0$, 
$\a_n=\|A_{n*}\|$, $\beta_n=\|b_{n*}\|$, and  $\gamma_{I,m}=\|A_{I,m*}\|$. For the first system in (\ref{Eq. principal2}) it follows from (\ref{Eq. principal en x}) that $z_n:=\|y_n\|$ satisfies the inequalities \begin{equation}\label{E aux para z x}
\frac{z_n}{c}\leq \beta_n+ r[n-p]_{|q|}  z_{n-p}+\sum_{j=1}^{n-1} \a_{n-j} z_j+\sum_{k=2}^n \sum_{*_k} \prod_{{1\leq l\leq N}\atop{1\leq j\leq i_l}} z_{n_{l,j}}\,\cdot  \gamma_{I,m}.
\end{equation} Note that we used that $\|y_{l,n_{l,j}}\|\leq \|y_{n_{l,j}}\|=z_{n_{l,j}}$. Now, we choose the sequence 
$$M_n=([n]^!_{|q|})^{1/p},$$ which is log-convex, see Remark \ref{Rmk. logcovex} and the inequality (\ref{Eq. Inq qbinom}). Recalling (\ref{Eq. Inq n-p qfactor}) we find that \begin{equation}\label{Eq. aux zxII}
\frac{z_n}{c\cdot M_n}\leq \frac{\beta_n}{M_n}+r \frac{z_{n-p}}{M_{n-p}}+\sum_{j=1}^{n-1} \frac{\a_{n-j}}{M_{n-j}} \frac{z_j}{M_j}+\sum_{k=2}^n \sum_{*_k} \prod_{{1\leq l\leq N}\atop{1\leq j\leq i_l}} \frac{z_{n_{l,j}}}{M_{n_{l,j}}}\,\cdot  \frac{\gamma_{I,m}}{M_m}.
\end{equation} 

For the second system we use the corresponding inequality to (\ref{Eq. principal en x}), to obtain that $z_n=\|y_n\|$ satisfies (\ref{E aux para z x}) with $[n-p]_{|q|}$ replaced by $|q|^{n-p}$. Then we divide by the log-convex sequence $$M_n=|q|^{n^2/2p},$$ and noticing that $|q|^{n-p}/M_n\leq 1/M_{n-p}$, we arrive again to (\ref{Eq. aux zxII}).

\textbf{Step 4: The majorized problem}.
Define recursively $w_n$ by $w_1=\|y_1\|$ and \begin{equation}\label{E aux para w}
\frac{w_n}{c}= \frac{\beta_n}{M_n} + r w_{n-p}+\sum_{j=1}^{n-1} \frac{\a_{n-j}}{M_{n-j}} w_j+\sum_{k=2}^n \sum_{*_k} \prod_{{1\leq l\leq N}\atop{1\leq j\leq i_l}} w_{n_{l,j}}   \frac{\gamma_{I,m}}{M_m}.
\end{equation} It follows by induction that $z_n/M_n\leq w_n$, for all $n\geq 1$. Recursion (\ref{E aux para w}) is equivalent to assert that   $\hat{w}(\tau):=\sum_{n\geq 1}w_n \tau^n$ satisfies  \begin{equation}\label{E para w}
\frac{1}{c}w(\tau) = \widetilde{b}(\tau)+r\tau^p \hat{w}(\tau)+\widetilde{A}(\tau)\hat{w}(\tau)+\sum_{|I|\geq2}\widetilde{A}_I(\tau)\hat{w}^{|I|},
\end{equation} where  $$\widetilde{b}(\tau)=\sum_{n=1}^\infty \frac{\beta_n}{M_n} \tau^n,\quad \widetilde{A}(\tau)=\sum_{m=1}^\infty \frac{\a_{n}}{M_n} \tau^n,\qquad \widetilde{A}_I(\tau)=\sum_{m=0}^\infty \frac{\gamma_{I,m}}{M_m} \tau^m,$$ define analytic (in this case, actually entire) functions of $\tau$. Consider the map $$H(\tau,w)=\widetilde{b}(\tau)+\left(r\tau^p+\widetilde{A}(\tau)-\frac{1}{c}\right)w+\sum_{k=2}^\infty \left(\sum_{|I|=k} \widetilde{A}_{I}(\tau)\right) w^k,$$ which is analytic at $(\tau,w)=(0,0)$, due to  estimates of the form \begin{equation}\label{Gr AI}
\|A_I\|\leq K\delta^{|I|},\qquad |I|\geq 2,\end{equation} for some constants $K,\delta>0$. They hold thanks to the analyticity of $F$ at the origin. Therefore, we can apply the Implicit Function Theorem to $H$ since $\frac{\d H}{\d w}(0,0)=-\frac{1}{c}\neq0$. Thus we find a unique analytic solution $\widetilde{w}(\tau)$ of $H(\tau,\widetilde{w}(\tau))=0$ and $\widetilde{w}(0)=0$. But $\hat{w}=\widetilde{w}$, since both are formal solutions, and there is only one such solution. In conclusion,  $\hat{w}\in\C\{\tau\}$. Thus $\hat{{y}}$ satisfies $$\|y_n\|\leq C L^n ([n]^!_{|q|})^{1/p},\qquad \|y_n\|\leq CL^n |q|^{n^2/2p},$$ respectively for each system, and for some constants $C,L>0$. Note that when we return to the original variable $\e$ before rank reduction, the only effect in the previous bounds is to change the radius $|\eta|<r$ by $|\e|<r^{1/\a}$. This conclude the case for the variable $x$ and $p>0$.

\

\textit{Solution as a power series in $x$. Case $p=0$}. Here the first equation in (\ref{Eq. principal2}) has the form $\e x d_q(y)(x,\e)=F(x,\e,y)$. Solving for $y$ as a power series in $x$, and assuming already that $y_0(\e)=0$,  the coefficient $y_n(\e)$ is now determined by 
\begin{equation}\label{Eq. principal en x 22}
 y_{n}(\e)=\left(\e [n]_qI_N-A_{0*}(\e)\right)^{-1}\left[b_{n*}(\e)+\sum_{j=1}^{n-1} A_{n-j*}(\e) y_{j}(\e)+\cdots\right],\quad n\geq 1.
\end{equation}

In this case, the reduction on the radius on $\e$ comes from $\left(\e [n]_q I_N-A_{0\ast}(\e)\right)^{-1}$. Indeed, write  $A_0=A(0,0)=A_{0\ast}(0)$ and choose $0<r<1/4\|A_{0}^{-1}\|$ such that $\|A_{0}-A_{0\ast}(\e)\|<1/4\|A_{0}^{-1}\|$, for all $|\e|\leq r$. Then, for each $n$, if $|\e|\leq r/|[n]_q|$, it follows that $$\|A_{0}^{-1}\left(\e [n]_q I_N+A_{0}-A_{0\ast}(\e)\right)\|\leq \frac{1}{2}<1,$$ and using the inequality 
$\|(I-B)^{-1}\|\leq\frac{1}{1-\|B\|}$, valid for $\|B\|<1$ (Neumann series), for any matrix norm, we find that  
\begin{align*}
 \| \left(\e [n]_q I_N-A_{0\ast}(\e)\right)^{-1}\|&= \| \left(I_N-A_0^{-1}\left(\e [n]_q I_N+A_{0}-A_{0\ast}(\e)\right)\right)^{-1}\cdot (-A_0)^{-1}\|\\
 & \leq 2\|A_0^{-1}\|=c.
 \end{align*}

If  $z_n=\sup_{|\e|\leq r/|[n]_q|} |y_n(\e)|$, $\a_n=\sup_{|\e|\leq r/|[n]_q|} \|A_{*n}(\e)\|$,  $\beta_n=\cdots$, we find \begin{equation}\label{E aux para z x II}
\frac{z_n}{c}\leq \beta_n+\sum_{j=1}^{n-1} \a_{n-j} z_j+\cdots.
\end{equation} It follows that $$\sup_{|\e|\leq r/|[n]_q|} |y_n(\e)|\leq CA^n,$$ for some constants $C,A>0$, and all $n$, as needed. Note we can interchange the terms $[n]_q$ and $q^n$ in this supremum by recalling (\ref{Eq. nqn asym}) and reducing $r$ if necessary. 

For the second equation in (\ref{Eq. principal2}), the coefficient $y_n(\e)$ is determined by (\ref{Eq. principal en x 22}) with $q^n$ instead of $[n]_q$. As before, the matrix $\left(\e q^n I_N-A_{0*}(\e)\right)^{-1}$ can be uniformly bounded for $|\e|\leq r/|q|^n$, for an adequate $r>0$. Letting $z_n=\sup_{|\e|\leq r/|q|^n} |y_n(\e)|$, $\a_n=\sup_{|\e|\leq r/|q|^n} \|A_{*n}(\e)\|$,  $\beta_n=\cdots$, we arrive again at (\ref{E aux para z x II}) and to the desired bounds.

To conclude, note that for the general case $\a>1$, when we return to the original variable $\e$, the condition $|\eta|\leq r/|q|^n$ means that $|\e|\leq r'/|q|^{\frac{n}{\a}}$ where $r'=r^{1/\a}$, since $\eta=\e^\a$.

\

\textit{Solution as a power series in $\epsilon$. Case $p\geq 0$}. Consider the system in (\ref{Eq. principal2}) with $\a=1$ and the search for the coefficients $u_n(x)$ of $\hat{y}$, as power series in $\e$. The coefficient $u_0(x)$ is determined by solving $F(x,0,u_0(x))=0$. After the change of variables $y\mapsto y-u_0(x)$ in equation (\ref{Eq. principal2}), we obtain a similar one with $F(x,0,0)=b_{*0}(x)=0$. As before, replacing $\hat{y}=(\hat{y}_1,\dots,\hat{y}_N)$,  $\hat{y}_j=\sum_{n=1}^\infty u_{j,n}(x) \e^n$, into the first system in (\ref{Eq. principal2}), we obtain the family of $q$-difference  equations \begin{align}\label{Eq. principal en e}
x^{p+1}d_q(u_{n-1})=&b_{*n}+\sum_{j=1}^n A_{*n-j}u_{j}
+\sum_{k=2}^n \sum_{\ast_k} \prod_{{1\leq l\leq N}\atop{1\leq j\leq i_l}} u_{l,n_{l,j}}\,\cdot  A_{*I,m},
\end{align} where the sum $(\ast_k)$ has the same structure as before. For the second system we obtain (\ref{Eq. principal en e}) with $x^{p}\sigma_q(u_{n-1})$ instead of $x^{p+1}d_q(u_{n-1})$.

Here $A(x,\e)=\sum_{n=0}^\infty A_{*n}(x)\e^n$ and similarly for the other terms. The radius $r>0$ is such that 
$A_{*n}\in \mathcal{O}_b(D_r)^{N\times N}$, $b_{*n}, A_{*I,n}\in \mathcal{O}_b(D_r)^N$, and $A_{*0}(x)$ is invertible and bounded for $|x|<r$. Thus, in both cases the coefficient $u_n$ are uniquely determined by the previous terms. Moreover, since the operators $d_q$ and $\sigma_q$ reduce the radius of convergence by a factor of $q$, $u_n\in \mathcal{O}_b(D_{r/|q|^n})$.

For Step $3$, let $c=\|A_{*0}^{-1}\|_0$, $z_n=\|u_{n}\|_n$, $\a_n=\|A_{*n}\|_n$,  $\beta_n=\|b_{*n}\|_n$, and $\gamma_{I,m}=\|A_{*I,m}\|_m$. Equation (\ref{Eq. principal en e}) and the properties of the $q$-Nagumo norms developed in Section \ref{Sec. Nagumo norms} establish that  
\begin{equation}\label{Eq. zn eI}
\frac{z_n}{c}\leq \beta_n+2 r^{p+1} z_{n-1}+\sum_{j=1}^{n-1} \a_{n-j}z_j+\sum_{k=2}^n \sum_{*_k} \prod_{{1\leq l\leq N}\atop{1\leq j\leq i_l}} z_{n_{l,j}}\,\cdot  \gamma_{I,m}.
\end{equation} In fact, the $q$-derivative $d_q$ is controlled by the inequality $$\|x^{p+1}d_q(u_{n-1})\|_n\leq \sup_{|x|\leq \frac{r}{|q|^n}} |x^{p+1}|\cdot \|d_q (u_{n-1})\|_n\leq 2r^{p+1} \|u_{n-1}\|_{n-1},$$ valid for $p\geq 0$. In the case of the operator $\sigma_q$ we have that $$\|x^{p}\sigma_q(u_{n-1})\|_n\leq \sup_{|x|\leq \frac{r}{|q|^n}} |x^{p}|\cdot \|\sigma_q (u_{n-1})\|_n\leq r^{p}\cdot r \|u_{n-1}\|_{n-1},$$ thus arriving to the same inequalities for the $z_n$.

Now we choose $M_n=1$ (no division at all) to conclude in the same way as in previous cases, that $$\|u_n\|_n\leq DB^n,$$ for some constants $B,D>0$. The definition of the $q$-Nagumo norms leads to $$|u_{n}(x)|\leq \frac{\|u_n\|_n}{d_n(|x|)^n}\leq \frac{DB^n }{(r-\rho)^n} = D L^n ,$$ where $L=B|q|/r(|q|-1)$ and  $|x|\leq {\rho}/{|q|^n}={r}/{|q|^{n+1}}$, as required, as we wanted to prove.

Finally, we  establish the case $\a>1$. If we return to the original variable $\e$ before rank reduction, $y(x,\e)=\sum_{j=0}^{\a-1} y_j(x,\epsilon^\a) \e^j=\sum_{n=0}^\infty u_n(x) \e^n$, and $y_j(x,\eta)=\sum_{k=0}^\infty u_{j,k}(x) \eta^k$, we have that $$u_{\a k+j}(x)=u_{j,k}(x),\qquad  j=0,1,\dots,\a-1, k\geq 0.$$ The proof has shown that $|u_{j,k}(x)|\leq D L^k$, for all $|x|\leq \rho/|q|^{k}$, and all $k$. Therefore, if $n\geq 0$ is divided by $\a$ and $n=\a k+j$, then $k=\lfloor n/\a\rfloor$. From here it is easy to conclude that $$\sup_{ |x|\leq {\rho}/{|q|^{\lfloor \frac{n}{\a}\rfloor}} } |u_n(x)|\leq C'L'^n,$$ for some $C'=C'(q)$, $L'=L'(q)>0$. In general, the value $r/|q|^{\lfloor n/\a\rfloor}$ means that the radius of the domain of the coefficients $u_n$ is reduced by a factor of $|q|$ every $\a$ steps, i.e., $u_{k\a}, u_{k\a+1},\dots, u_{k\a+\a-1}\in\mathcal{O}(D_{r/|q|^k})$, for all $k$.
\end{proof}

Finally, we conclude by proving Theorem \ref{Prop1}.

\begin{proof}[Proof of Theorem \ref{Prop1}] We proceed as in the proof of Theorem \ref{Thm. Main} searching for a solution as a power series in $\e$. In this case, after rank reduction, i.e., assuming $\a=1$, if we search for a solution of (\ref{Eq.3}) of the form $\hat{y}=\sum_{n=0}^\infty u_n(x) \e^n$, we will arrive at recurrence (\ref{Eq. principal en e}) with $p=-1$. Now, the induced inequalities for $z_n=\|u_n\|_n$ are 
\begin{equation*}
\frac{z_n}{c}\leq \beta_n+2 |q|^n z_{n-1}+\sum_{j=1}^{n-1} \a_{n-j}z_j+\cdots.
\end{equation*} In this case we choose $$M_n=|q|^n\cdot |q|^{n^2/2},$$ which  is log-convex since $M_n/M_{n-1}=|q|^{n+\frac{1}{2}}$ is increasing in $n$, and satisfies $|q|^n/M_n\leq 1/M_{n-1}$. Therefore, dividing the previous inequality by $M_n$ and proceeding as before, we find that $$\| u_n\|_n\leq DB^n |q|^n |q|^{n^2/2},$$ for some constants $B,D>0$. Recalling the definition of $q$-Nagumo norms, this means that $$|u_n(x)|\leq DL^n |q|^{n^2/2},$$ for $|x|\leq \rho/|q|^{n}=r/|q|^{n+1}$, where $L=B|q|^2/r(|q|-1)>0$.

Finally, for general $\a>1$, using the same notation that in the last paragraph of the previous proof, in this case $\sup_{|x|\leq \rho/|q|^k} |u_{j,k}(x)|\leq DL^k |q|^{k^2/2}$. Therefore, if $n=\a k+j$, and $k=\lfloor n/\a \rfloor$, we find that $$\sup_{|x|\leq \rho/|q|^{\lfloor n/\a \rfloor}} |u_n(x)|\leq DL^{n/\a} |q|^{(n/\a)^2/2},$$ as required. Here  we have assumed that $L>1$. This concludes the proof.
\end{proof}

\begin{nota}\label{Rmk RR} If we reduce rank in $x$ and $\epsilon$, we should consider the decomposition $$y(x,\e)=\sum_{{0\leq l<p}\atop {0\leq j< \a}} y_{l,j}(z,\eta) x^l\e^j,\qquad \text{ where } z=x^p, \eta=\e^\a.$$ Computing $d_{q,x}(y)$, it follows using (\ref{qLeib}) and (\ref{Eq. npq}) that $$\epsilon^{\a} x^{p+1}d_{q,x}(y)=\sum_{{0\leq l<p}\atop {0\leq j< \a}} \left([p]_q q^l \eta z^2 d_{q^p,z}(y_{l,j})(z,\eta) +[l]_q \eta z y_{l,j}(z,\eta)\right) x^l\e^j.$$ Then $w(z,\eta)=(y_{l,j})_{0\leq l<p, 0\leq j<\a}$, written in lexicographical order, satisfies  $$\eta z^2 d_{q^p,z}w=G(z,\eta,w),$$ where $D_w(0,0,0)=[p]_q^{-1}D_q\widetilde{A}$, $D_q=\text{diag}(I_{\a N}, q^{-1} I_{\a N},\dots,q^{-(p-1)} I_{\a N})$, and $\widetilde{A}$ is a lower-triangular matrix having $A(0,0)$ as diagonal blocks. Apart from the linear part, $G$ also contains a term depending on $q$, namely, $\eta z[p]_q^{-1}D_qM_q$, where $$M_q=\text{diag}(0, [1]_q I_{\a N},\dots,[p-1]_q I_{\a N}).$$ Therefore, this reduction would force to change the nature of the problem by introducing $q$ on the non-$q$-difference part of the equation.
\end{nota}

\section{Examples}\label{Sec. Examples}

This section is devoted to give examples of our results. We remark that in few cases the formal power series solutions can be easily computed.

%{\color{red}{$$\e x d_q(y)(x,\e)=y-\frac{x}{1-x}$$

%$$\hat{y}=\sum_{m=1}^\infty \frac{x^m}{1-[m]_q\e}=\sum_{m\geq1, n\geq 0} [m]_q^n x^m \e^n$$

%If $q\to 1$, $\e x \d_x y=y-\frac{x}{1-x}$, $\hat{y}=\sum_{m=1}^\infty \frac{x^m}{1-m\e}=\sum_{m\geq1, n\geq 0} m x^m \e^n$. Solution to homogeneous: $x^{1/\e}=e^{\log(x)/\e}$.

%What is the $q$-analogue to $\e xd_qy=y$?}}

\begin{eje} ($\sigma_q$, case $p=0$) Consider the problem $$\e^\a \sigma_{q,x}y=A(x)y-b(x),$$ where $\a\in\N^+$,  $y\in\C^N$ and $A(0)$ is an invertible matrix. It follows that $$\hat{y}_q(x,\e)=\sum_{m=0}^\infty \left[A(q^m x)\cdots A(qx)A(x)\right]^{-1} b(q^m x) \e^{\a m},$$ is the unique formal power series solution of the problem. We see that $\hat{y}_q\in\mathcal{O}_{0}^{q^{1/\a}}$,  confirming Theorem \ref{Thm. Main}(2) for $\sigma_q$. In the limit $q\to 1$, the series reduces to $\hat{y}_1=\left[A(x)-\e^\a I_N\right]^{-1} \cdot b(x)$, which is of course the unique analytic solution of the limit problem $\e^\a y(x,\e)=A(x)y(x,\e)-b(x)$. A particular interesting case is the equation $$\e^\a \sigma_{q,x}y=(1-x)y-1,$$ generating the solution  \begin{align*}
%\sum_{m=0}^\infty \frac{\e^{\a n}}{(1-x)(1-qx)\cdots(1-q^m x)}
\hat{y}_q(x,\epsilon)=\sum_{m=0}^\infty \frac{\e^{\a m}}{(x;q)_{m+1}}=\sum_{n,m=0}^\infty  \qbinom{n+m}{m} x^n \e^{\a m}=\sum_{n=0}^\infty \frac{x^n}{(\e^\a;q)_{n+1}}.
\end{align*} Note we used Heine's binomial formula, see \cite[p. 28]{Kac}. Since $\qbinom{n+m}{m}$ is a monic polynomial in $q$ of degree $mn$, we have the precise bounds given in Theorem \ref{Thm. Main}(2) for $\sigma_q$. Additionally, the series reduces to $\hat{y}_1=(1-x-\e^\a)^{-1}$ for $q=1$.

We also highlight the case $A(x)=I_N$ and $b(x)=\sum_{m=0}^\infty a_m x^m\in\C\{x\}$ for which $$\hat{y}_q(x,\e)=\sum_{m=0}^\infty b(q^m x) \e^{\a m}=\sum_{n,m=0}^\infty a_n q^{nm} x^n\e^{
\a m}=\sum_{n=0}^\infty \frac{a_n x^n}{1-q^n \e^\a}.$$ If $\a=1$ and  $b(x)=(1-x)^{-1}$, we recover the series $M(x,\e)$ of Example \ref{Ex. Model}.
\end{eje}

\begin{eje}($d_q$, case $p=0$)
Consider the problem 
$$\e^\a x d_{q,x}y=y-f_0(x),$$
with $f_0(x)=\sum_{n=0}^\infty a_n x^n\in\C\{x\}$. Solving for $\hat{y}=\sum_{m=0}^{\infty} u_m(x)\e^{\a m}$ we see that $u_m(x)=xd_q(u_{m-1})=\cdots=(xd_q)^m(f_0)$. Therefore, $$  \hat{y}=\sum_{n\geq 1,m\geq 0} a_n [n]_q^m x^n \e^{\a m}=\sum_{n=1}^\infty \frac{a_n x^n}{1-[n]_q \e^\a} \in \mathcal{O}_0^{q^{1/\a}},$$ as claimed by Theorem \ref{Thm. Main}(2). The same conclusion can be achieved if we have solved the problem writing $\hat{y}=\sum_{n=0}^\infty y_n(\e) x^n$, by solving the recurrence $[n]_q \e^\a y_n(\e)=y_n(\e)-a_n$.

% $$u_n(x)=\frac{1}{(q-1)^m}\sum_{j=0}^{m}\binom{m}{j}(-1)^ju_{n-m}(q^{m-j}x)==\frac{1}{(q-1)^n}\sum_{j=0}^{n}\binom{n}{j}(-1)^ju_{0}(q^{n-j}x),\quad 0\le m\le n.$$ Assume that $r>0$ is the radius of convergence of the power series defining $u_0$ at the origin. Then, there exists $C>0$ such that $$|u_0(q^{n-j}x)|\le \frac{C}{1-\frac{q^{n-j}|x|}{r}},$$ for every fixed $0\le j\le n$ and $0\le |x| < r/q^{n-j}$. Therefore, $u_n$ is holomorphic for $|x|<r/q^n$, which agrees with Theorem~\ref{Thm. Main}, for the situation $\alpha=1$ and $p=0$. In addition to this, for every $x$ with $|x|<r/q^n$ one has that $$|u_{n}(x)|\le \frac{Cr}{(q-1)^n}\sum_{j=0}^{n}\binom{n}{j}\frac{1}{r-q^{n-j}r/q^{n}}\le \frac{C}{(q-1)^n}2^n\frac{q^n}{q^{n}-r}.$$ On the other hand, the formal solution of the equation can be written in the form $\hat{y}=\sum_{n=0}^{\infty}y_n(\epsilon)x^n$. Let us write $u_0(x)=\sum_{n=0}^{\infty}y_{0n}x^n$. One has $y_0(\epsilon)=y_{00}$ and for every $n\ge1$ $$\epsilon y_n(\epsilon)[n]_{q}=y_{n}(\epsilon)-y_{0n}.$$ Therefore, $y_n(\epsilon)=\frac{y_{0n}}{1-\epsilon [n]_q}$ is holomorphic on $0\le |\epsilon|<\frac{1}{[n]_q}\le\frac{C(q-1)}{q^n}$ (see~\ref{Eq. nqn asym}) for some $C>0$. For such $\epsilon$, analogous arguments determine positive constants $A,B$ such that $|y_n(\epsilon)|\le AB^n$ for all $|\epsilon|\le r_0/q^n$, for some fixed $r_0>0$. The estimate provided in Theorem~\ref{Thm. Main} is attained.
\end{eje}

\begin{eje} ($\sigma_q$, $p>0$, $\a=1$) Consider the equation $$\epsilon x^p\sigma_{q,x}y=y-f(x,\epsilon),$$ where $f(x,\epsilon)=\sum_{m=0}^\infty f_m(x) \e^m\in \C\{x,\e\}$. It follows that $\hat{y}=\sum_{m=0}^\infty u_m(x) \e^m$ satisfies the initial problem if and only if $u_0(x)=f_0(x)$ and $$x^p u_{n-1}(qx)=u_n(x)-f_n(x),\qquad n\geq 1.$$ We can solve this recursively to find that $$u_n(x)=\sum_{j=0}^n q^{pj(j-1)/2} x^{jp} f_{n-j}(q^j x).$$ These coefficients exhibit the growth $$\sup_{|x|\leq r/|q|^n} |u_n(x)|\leq CA^n,$$ for some fixed $r>0$, due to the restriction on the domain of $x$. This holds in particular for the case $f(x,\epsilon)=f_0(x)$ for which $u_n(x)=q^{pn(n-1)/2} x^{pn} f_0(q^n x)$. For instance, the equation $$\e x^2\sigma_{q,x}(y)=y-x$$ has the unique formal  solution $$\hat{y}(x,\e)=\sum_{n=0}^\infty q^{n(n-1)} x^{2n}(q^nx)\e^n=\sum_{n=0}^\infty (q^{n^2}\e^n) x^{2n+1}.$$ Therefore, the coefficients $y_k(x)=y_{2n+1}(\e)=q^{n^2}\e^n$ grow as $$\sup_{|\e|\leq r} |y_{k}(\e)|= r^n|q|^{n^2}=r^{\frac{k-1}{2}}|q|^{(k-1)^2/2\cdot 2}.$$ This confirms the bounds found in Theorem \ref{Thm. Main}(1) for $\sigma_q$, showing that they are optimal. In terms of the scalar coefficients of $\hat{y}=\sum a_{n,m} x^n \e ^m$, $a_{2n+1,n}=q^{n^2}$ and $a_{n,m}=0$ otherwise. In this case, $$\min\{|q|^{(2n+1)^2/2},|q|^{(2n+1)n}\}=|q|^{(2n+1)^2/2}=|q|^{n^2+n+\frac{1}{4}},$$ thus $\hat{y}$ belongs precisely to $\mathcal{O}^q_{x^2\e}$.
\end{eje}

\begin{eje} ($d_q$, $p=\a=1$) Consider the equation $$\epsilon x^2d_{q}y=y-f_0(x),$$ where $f_0(x)=\sum_{n=0}^\infty a_n x^n \in \C\{x\}$. Searching  for a solution of the form  $\hat{y}=\sum_{m=0}^{\infty} u_m(x)\e^{m}$ we see that $u_m(x)=x^2d_q(u_{m-1})=\cdots=(x^2d_q)^m(f_0)$. Therefore, $$  \hat{y}=\sum_{n\geq 1,m\geq 0} a_n [n]_q[n+1]_q\cdots [n+m-1]_q x^{n+m} \e^{m}=\sum_{n\geq 1, m\geq 0} a_{n,m} x^{n+m} \e^m,$$ since $(x^2d_q)^m(x^n)=[n]_q[n+1]_q\cdots [n+m-1]_q x^{n+m}$. Recalling (\ref{Eq. nqn asym}) we see that $|a_{n+m,m}|\sim q^{nm}\cdot q^{m(m-1)^2/2}/(q-1)^{m}$ for $n,m\to+\infty$. Therefore, up to a geometric terms (depending on $q$), we see that $$|a_{n+m,m}|\sim q^{nm}\cdot q^{m^2/2}.$$ On the other hand, Theorem \ref{Thm. Main}(1) provides bounds of type $$\min\{|q|^{(n+m)^2/2}, |q|^{(n+m)m}\}=\begin{cases}
|q|^{m(n+m)}, & m\leq n,\\
|q|^{(n+m)^2/2}, & n\leq m. 
\end{cases}$$ In particular, if $n=m$ the theorem asserts bounds of the form $|q|^{2m^2}$ but the actual term grows as $|q|^{3m^2/2}$, which is smaller. We will see in the following example that the bounds provided by our results are attained in this case as well.
\end{eje}

\begin{eje} ($d_q$, $p=\a=1$) Consider the scalar equation \begin{equation}
\e  x^2 d_q(y)=(1+x)y - x\e.  
\end{equation} It has the unique formal power series solution $\hat{y}(x,\e)\in\C[[x,\e]]$ given by $$\hat{y}(x,\e)=\sum_{n=1}^\infty y_n(\e) x^n=\sum_{1\leq m\leq n} a_{n,m} x^n \e^m,\qquad y_n(x)=\e\ \prod_{j=1}^{n-1} ([j]_q\e-1).$$ Therefore, $a_{n,n}=[n-1]^!_q$ and the bounds for the coefficients are attained. Indeed, $\min\{|q|^{n^2/2},|q|^{nm}\}=|q|^{n^2/2}$ for $n=m$. On the other hand, the solution is not easy to write when we expand $\hat{y}=\sum_{m=1}^\infty u_m(x) \e^m$. In this case the recurrences are  $$u_1(x)=\frac{x}{1+x},\qquad u_m(x)=\frac{x^2}{1+x} d_q(u_{m-1})(x),\quad m\geq 2.$$ Therefore, $$u_m(x)=\frac{x^m\cdot P_m(x,q)}{(1+x)^m(1+qx)^{m-1}(1+q^2x)^{m-2}\cdots (1+q^{m-1}x)},$$ where $P_m(x,q)\in \C[x,q]$. They can be found using the recursion \begin{align*}
P_{m+1}(x,q)=&\frac{q^m(1+x)^m P_m(qx,q)-(1+qx)(1+q^2x)\cdots(1+q^mx)P_m(x,q)}{q-1}\\
= & q^m(1+x)^m d_qP_m(x,q)+ \frac{q^m(1+x)^m-(-qx;q)_m}{q-1}P_m(x,q).
\end{align*} After some computer calculations, the first values of $P_m$ are found to be $P_2(x,q)=1$, 
$P_3(x,q)=-q^{2} x^{2} + q x + q + 1$, and 

\begin{align*}
P_4(x,q)& =   q^{7} x^{5} + x^{4} \left(- 3 q^{6} - 2 q^{5}\right) + x^{3} \left(- 4 q^{6} - 6 q^{5} - 3 q^{4} - 2 q^{3}\right) \\
&+ x^{2} \left(- q^{6} - 2 q^{5} - q^{4} - q^{3}\right) + x \left(q^{4} + 3 q^{3} + 4 q^{2} + 2 q\right) + + q^{3} + 2 q^{2} + 2 q + 1,\\
P_5(x,q)& =  - q^{16} x^{9} - x^{8} \left(- 6 q^{15} - 7 q^{14} - 3 q^{13}\right)\\
& - x^{7} \left(- 10 q^{15} - 19 q^{14} - 14 q^{13} - 8 q^{12} - 5 q^{11} - 3 q^{10}\right)\\
&- x^{6} \left(- 5 q^{15} - 11 q^{14} - 4 q^{13} + 7 q^{12} + 12 q^{11} + 8 q^{10} + 6 q^{9} + q^{8}\right)\\
& - x^{5} \left(- q^{15} - 2 q^{14} + 8 q^{13} + 35 q^{12} + 64 q^{11} + 71 q^{10} + 61 q^{9} + 40 q^{8} + 19 q^{7} + 6 q^{6}\right)\\
& - x^{4} \left(3 q^{13} + 22 q^{12} + 55 q^{11} + 84 q^{10} + 98 q^{9} + 93 q^{8} + 69 q^{7} + 37 q^{6} + 12 q^{5} + 3 q^{4}\right)\\
& - x^{3} \left(4 q^{12} + 15 q^{11} + 30 q^{10} + 44 q^{9} + 54 q^{8} + 50 q^{7} + 34 q^{6} + 18 q^{5} + 8 q^{4} + 2 q^{3}\right)\\
& - x^{2} \left(q^{11} + 3 q^{10} + 5 q^{9} + 5 q^{8} - 10 q^{6} - 15 q^{5} - 13 q^{4} - 8 q^{3} - 2 q^{2}\right)\\
&- x \left(- q^{8} - 4 q^{7} - 11 q^{6} - 18 q^{5} - 21 q^{4} - 18 q^{3} - 10 q^{2} - 3 q\right)\\
& + q^{6} + 3 q^{5} + 5 q^{4} + 6 q^{3} + 5 q^{2} + 3 q +1,
\end{align*} The following term, $P_6(x,q)$ has $203$ terms and leading term $x^{14}q^{30}$. It can be shown that 
$$\text{deg}_x(P_m)=\frac{m(m-1)}{2}-1,\qquad \text{deg}_q(P_m)=\frac{(m-1)(m-2)(m+3)}{6},$$ but no easy closed formula seems to be available in the general case.

\end{eje}

We conclude this section with examples of Theorem \ref{Prop1}.

\begin{eje} Consider the scalar equation \begin{equation}\label{Eq. Ex6}
   \e^\a d_{q,x}y=y-\frac{1}{1-x}. 
\end{equation} Letting $\eta=\e^\a$ and searching for a formal solution of the form $\hat{y}_q=\sum_{n=0}^\infty u_n(x)\eta^n=\sum_{n=0}^\infty v_{\a n}(x)\e^{\a n}$, we find that $u_0(x)=\frac{1}{1-x}$, and $$u_{n+1}(x)=d_{q,x}u_n(x),\qquad n\geq 0.$$ It follows by induction that $$u_n(x)=\frac{[n]^!_q}{(x;q)_{n+1}}.$$ Thus, for every $0<r<1$, $u_n\in\mathcal{O}_b(D_{r/|q|^n})$ and $$\sup_{|x|\leq r/|q|^n} |u_n(x)|\leq  \frac{|[n]^!_q|}{(1-r)^{n+1}}.$$ Therefore, $\hat{y}$ has the precise growth described in Theorem \ref{Prop1}, since the previous bounds can be written as $$\sup_{|x|\leq r/|q|^{\a n/\a}} |v_{\a n}(x)|\leq CA^{\a n} |q|^{(\a n)^2/2\a ^2},$$ for adequate constants $C,A>0$. Also note that the expansion of $\hat{y}_q$ is $$\hat{y}_q=\sum_{m=0}^\infty \frac{[m]_q^!}{(x;q)_{m+1}} \e^{\a m}=\sum_{n,m=0}^\infty \frac{[n+m]_q^!}{[n]_q^!} x^n \e^{\a m}.$$ Thus the coefficients satisfy $a_{n,\a m}=[n+m]_q^! [n]_q^!\sim |q|^{(n+m)^2/2}/|q|^{n^2/2}=|q|^{nm+m^2/2} $ as expected.
\end{eje}

\begin{eje}
More generally, consider the problem (\ref{Eq. Ex6}) with $f\in\mathcal{O}(D_r)$ instead of $1/(1-x)\in\mathcal{O}(D_1)$. We obtain the formal solution $$\hat{y}(x,\epsilon)=\sum_{n=0}^\infty d_q^n f(x) \epsilon^n.$$ In this case, the Borel transform in $\e$ of $\hat{y}$ is $$\mathcal{B}_q(\hat{y})(\xi,x):=\sum_{n=0}^\infty \frac{d_q^nf(x)}{[n]_q^!} \xi^n.$$ In the differential case, i.e., $q=1$, we see that $\mathcal{B}_1(\hat{y})(\xi,x)=f(x+\xi)$, $|x+\xi|<r$ thanks to Taylor's theorem. However, in the $q$-difference case, in general only for $x=0$, we have that $\mathcal{B}(\hat{y})(\xi,0)=f(\xi)$, $|\xi|<r$. For instance, if we take the $q$-exponential $f(x)=e_q(x):=\sum_{n=0}^\infty x^n/[n]_q^!$, we find that $$\mathcal{B}_q(\hat{y})(\xi,x)=e_q(x)e_q(\xi)\neq e_q(x+\xi),$$ if $x\neq0, \xi\neq 0$, see \cite{Kac}. Note we used that $e_q$ is the unique solution of the problem $d_qf=f$ and $f(0)=1$,
\end{eje}

\section{Appendix: \texorpdfstring{$q$}--Nagumo norms suitable for confluence}\label{Sec. Nagumo norms I}

In this section we include non-modified $q$-Nagumo norms suitable for confluence. If fact, for $q=1$ they reduce to the usual non-modified Nagumo norms. These are useful to recover the known results for the rate of divergence of equation (\ref{Eq. MS}).

Let $q>1$ and $f\in\mathcal{O}_b(D_r)$. The $n$th $q$-Nagumo norm is defined by $$\|f\|'_n:=\sup_{|x|\leq r/q^n} |f(x)| (r-q^n|x|)^n.$$ Note that if $\|f\|'_n$ is finite, we have that $$|f(x)|\leq \frac{\|f\|'_n}{(r-q^n|x|)^n}\leq \frac{\|f\|'_n}{(r-\rho)^n},$$ for $|x|\leq \rho/q^n<r/q^n$. For these norms we have the following properties.

\begin{lema} Let $q>1$, $n,m\in\N$ and $f,g \in\mathcal{O}_b(D_r)$. Then:\begin{enumerate}
    \item $\|f+g\|'_n\leq \|f\|'_n+\|g\|'_n$ and $\|fg\|'_{n+m}\leq \|f\|'_n \|g\|'_m$.
    
    \item $\|d_q(f)\|'_{n+1}\leq eq^n(n+1)\|f\|'_n$.
    
    \item $\|\sigma_q(f)\|'_{n+1}\leq r\|f\|'_n$.
\end{enumerate}
\end{lema}

\begin{proof} For (1) assume $n\leq m$. Since $r-q^{n+m}|x|\leq r-q^m|x|\leq r-q^n|x|$ and $r/q^{n+m}\leq r/q^m\leq r/q^n$, we find $$|f(x)g(x)|(r-q^{n+m}|x|)^{n+m}\leq |f(x)||g(x)|(r-q^n|x|)^n(r-q^m|x|)^m\leq\|f\|'_n \|g\|'_m,$$ for all $|x|\leq r/q^{n+m}$. This proves the inequality for the product. 

For (2) fix $x$ such that $|x|\leq r/q^{n+1}$ and choose $\rho>0$ with $q|x|+\rho<r/q^n$. Notice that if $1\leq t\leq q$, then $|tx|+\rho<r/q^n$. Therefore, the disk $D(tx,\rho)$ is contained in $D(0,r/q^n)$. By Cauchy's formulas we can write \begin{align*}
    |d_q(f)(x)|=\left|\frac{1}{(q-1)} \int_1^q f'(tx)dt \right|=\left|\frac{1}{2\pi i(q-1)} \int_1^q \int_{\gamma} \frac{f(\xi)}{(\xi-tx)^2}d\xi dt \right|, 
\end{align*} where $\gamma$ is the boundary of $D(tx,\rho)$. It follows that $$|f'(tx)|\leq \frac{1}{\rho}\sup_{|\xi-tx|=\rho} |f(\xi)|\leq \frac{{\|f\|'_n}/{\rho}}{(r-q^n(t|x|+\rho))^n}\leq \frac{{\|f\|'_n}/{\rho}}{(r-q^n(q|x|+\rho))^n}.$$ Therefore, $$|d_q(f)(x)|\leq \frac{{\|f\|'_n}/{\rho}}{(r-q^{n+1}|x|-q^n\rho)^n}.$$ Choosing $q^n\rho=\frac{r-q^{n+1}|x|}{n+1}$ we find that $$(r-q^{n+1}|x|)^{n+1}|d_q(f)(x)|\leq \frac{q^n(n+1)\|f\|'_n}{(1-\frac{1}{n+1})^n},$$ and the result follows from the well-known inequality $(1-\frac{1}{n+1})^{-n}<e$.

To prove (3) simply note that by definition of $\|f\|'_n$ we have $$|f(qx)|\leq \frac{\|f\|'_n}{(r-q^{n}|qx|)^{n}},\qquad \text{ for } |qx|<r/q^{n}.$$ Therefore, $|f(qx)|(r-q^{n+1}|x|)^{n+1}\leq \|f\|'_n (r-q^{n+1}|x|)<r\|f\|'_n$ from which we conclude the desired inequality.
\end{proof}

\begin{coro}[Proposition 5 \cite{CDMS}] Consider the system of doubly-singular differential equations (\ref{Eq. MS}) with the same hypotheses on $F$ as in Theorem \ref{Thm. Main}. Then the problem has a unique formal power series solution $\hat{y}=\sum_{n,m\geq0} a_{n,m} x^n \e^m\in\C[[x,\e]]^N$ which is $1$-Gevrey in the monomial $x^p\e^\sigma$, i.e., there are constants $C,A>0$ such that \begin{equation}\label{Eq.Last}
    |a_{n,m}|\leq CA^{n+m} \min\{n!^{1/p}, m!^{1/\a}\}.
\end{equation}
\end{coro}

\begin{proof} For each $q>1$  consider the equation (\ref{Eq. principal}) for $\a=1$. We know by Theorem \ref{Thm. Main} that there is a unique formal power series solution $\hat{y}_q=\sum_{n=0}^\infty y_{n}^{[q]}(\e) x^n=\sum_{m=0}^\infty u_{m}^{[q]}(x) \e^m=\sum_{n,m\geq 0} a_{n,m}(q) x^n\e^m \in\C[[x,\e]]^N$. By letting $q\to 1^+$, it is easy to check that $\hat{y}=\lim_{q\to 1^+} \hat{y}_q$ (term by term) as in (\ref{Eq. y ini}) satisfies (\ref{Eq. MS}).

A second look at the proof of Theorem \ref{Thm. Main}, Steps 3 and 4, shows that if we modify $\widetilde{b}, \widetilde{A}$, and each $\widetilde{A}_I$ changing $M_n=([n]!_{q})^{1/p}$ by $1$, the same argument works to conclude that $$\sup_{|x|\leq r}|y_n^{[q]}(\e)|\leq CL^n ([n]!_{q})^{1/p},$$ for constants $C,L>0$ \textit{independent} of $q$. Therefore, letting $q\to 1^+$, we find that $$\sup_{|x|\leq r} |y_n(\e)|\leq CL^n n!^{1/p}.$$

To obtain the bounds as a power series in $\e$, we repeat the corresponding proof but now using $z_n=\|u_n^{[q]}\|_n'$. In this case, the inequalities (\ref{Eq. zn eI}) take the form $$\frac{z_n}{c}\leq \beta_n+er^{p+1} n z_{n-1}+\cdots.$$ Thus, after dividing by $M_n=n!$ and arguing as usual, we conclude that $\|u_n^{[q]}\|_n\leq DB^n n!$, for some constants $D,B>0$ independent of $q$. Therefore, $$\sup_{|\e|\leq \rho/q^n} |u_n^{[q]}(x)|\leq \frac{DB^n}{(r-\rho)^n} n!,\quad \text{ for a fixed } 0<\rho<r.$$ Letting $q\to 1^+$, we find that $\sup_{|x|\leq \rho} |u_n(x)|\leq \frac{DB^n}{(r-\rho)^n} n!$, as required. Undoing the rank reduction by using the bounds $$|u_n(x)|=|u_{\a k +j}(x)|\leq D'T^k k!\leq D' T^{n/\a} (\a k)!^{1/\a}\leq D' T^{n/\a} n!^{1/\a},$$ we obtain the general case. Finally, Cauchy's inequalities lead to (\ref{Eq.Last}).
\end{proof}

\begin{nota} If $n=\a k+j$, $0\leq j<\a$, it follows that  $$[k]!_{|q|}^\a\leq [\a k]^!_{|q|}\leq [n]^!_{|q|},$$ for $|q|>1$. Letting $q\to 1^+$ we can recover the usual inequality $k!\leq (\a k)!^{1/\a}\leq n!$ employed in the previous proof.
\end{nota}

\end{document}